\documentclass[12pt]{amsart}
\usepackage{amssymb,amscd}

\newtheorem{thm}[equation]{Theorem}
\newtheorem{lem}[equation]{Lemma}
\newtheorem{prop}[equation]{Proposition}
\newtheorem{cor}[equation]{Corollary}
\theoremstyle{definition}

\newtheorem{example}[equation]{Example}

\theoremstyle{remark}
\newtheorem{remark}[equation]{Remark}
\numberwithin{equation}{section}

\begin{document}
\title{Presenting Schur algebras}
\author{Stephen Doty}
\address{Loyola University Chicago, Chicago, Illinois 60626 U.S.A.}
\email{doty@math.luc.edu, tonyg@math.luc.edu}
\author{Anthony Giaquinto}

\date{March 31, 2002}

\allowdisplaybreaks

\begin{abstract} 
Motivated by work of R.M.\ Green, we obtain a presentation of Schur
algebras (both the classical and quantized versions) in terms of
generators and relations. The presentation is compatible with the
usual presentation of the (quantized or classical) enveloping algebra
of $\mathfrak{gl}_n$.  As a result, we obtain a new ``integral'' basis
for Schur algebras which is a subset of Kostant's basis of the
integral form of the enveloping algebra (or its $q$-analogue). 
Projection onto an appropriate component gives a new
``integral'' basis and a presentation for the Hecke algebra,
compatible with the basis and presentation for the Schur algebra.
Finally, we find a second presentation of Schur algebras which is
similar to Luzstig's modified form of the quantized enveloping
algebra.
\end{abstract}

\maketitle
\newcommand{\N}{{\mathbb N}}
\newcommand{\Z}{{\mathbb Z}}
\newcommand{\Q}{{\mathbb Q}}
\newcommand{\R}{{\mathbb R}}
\newcommand{\C}{{\mathbb C}}
\newcommand{\g}{{\mathfrak g}}
\newcommand{\n}{{\mathfrak n}}
\newcommand{\dist}{\operatorname{Dist}}
\newcommand{\per}{\operatorname{per}}
\newcommand{\cov}{\operatorname{cov}}
\newcommand{\non}{\operatorname{non}}
\newcommand{\cf}{\operatorname{cf}}
\newcommand{\add}{\operatorname{add}}
\newcommand{\End}{\operatorname{End}}
\newcommand{\Ext}{\operatorname{Ext}}
\newcommand{\Hom}{\operatorname{Hom}}
\newcommand{\Tor}{\operatorname{Tor}}
\newcommand{\ch}{\operatorname{ch}}
\newcommand{\ind}{\operatorname{ind}}
\newcommand{\coind}{\operatorname{Coind}}
\newcommand{\res}{\operatorname{res}}
\newcommand{\soc}{\operatorname{soc}}
\newcommand{\rad}{\operatorname{rad}}
\newcommand{\Aut}{\operatorname{Aut}}
\newcommand{\Dist}{\operatorname{Dist}}
\newcommand{\Lie}{\operatorname{Lie}}
\renewcommand{\ker}{\operatorname{Ker}}
\newcommand{\im}{\operatorname{im}}
\newcommand{\GL}{{\sf GL}}
\newcommand{\SL}{{\sf SL}}
\newcommand{\gl}{\mathfrak{gl}}
\renewcommand{\sl}{\mathfrak{sl}}
\newcommand{\B}{T}
\newcommand{\KK}{\mathcal{K}}
\newcommand{\divided}[2]{#1^{(#2)}}
\newcommand{\ct}{\chi}
\newcommand{\ctl}{\ct_L}
\newcommand{\ctr}{\ct_R}
\newcommand{\sqbinom}[2]{\begin{bmatrix}#1\\#2\end{bmatrix}}
\newcommand{\U}{\mathbf{U}}
\let\sect=\S
\renewcommand{\S}{\mathbf{S}}
\newcommand{\A}{\mathcal{A}}
\newcommand{\V}{\mathbf{V}}
\newcommand{\tgen}{t}
\newcommand{\ep}{\varepsilon}
\newcommand{\al}{\alpha}
\newcommand{\be}{\beta}
\newcommand{\Ea}{E_{\alpha}}
\newcommand{\Eb}{E_{\beta}}
\newcommand{\Eg}{E_{\gamma}}
\newcommand{\Ep}{E_{\rho}}
\newcommand{\Eij}{E_{ij}}
\newcommand{\Ekl}{E_{kl}}
\newcommand{\vbinom}[2]{\begin{bmatrix}#1\\#2\end{bmatrix}}
\newcommand{\vbinomm}[3]{\begin{bmatrix}#1; #2\\#3\end{bmatrix}}
\newcommand{\1}{1_{\lambda}}
\newcommand{\Lnd}{\Lambda (n,d)}
\newcommand{\Ua}{\mathbf{U}_{\mathcal{A}}}
\newcommand{\basis}{Y}
\newcommand{\qbasis}{\mathbf{Y}}

\parskip=2pt

\section*{Introduction} \label{sec:intro}

Let $R$ be a commutative ring.  The (classical) Schur algebra
$S_R(n,d)$ may be defined as the algebra $\End_{\Sigma_d}(V_R^{\otimes
d})$ of linear endomorphisms on the $d$th tensor power of an
$n$-dimensional free $R$-module $V_R$ commuting with the action of the
symmetric group $\Sigma_d$, acting by permutation of the tensor places
(see \cite{Green}).  The Schur algebras form an important class of
quasi-hereditary algebras, and, when $R$ is an infinite field, the
family of Schur algebras $\{ S_R(n,d) \}_{d\ge 0}$ determines the
polynomial representation theory of the general linear group
$\GL(V_R)$.

All these algebras, for various $R$, can be constructed from the
integral form $S_\Z(n,d)$ by base change, since $S_R(n,d) \cong
R\otimes_\Z S_\Z(n,d)$.  Fixing our base field at $\Q$ (we could use
any field of characteristic zero), we henceforth write $S(n,d)$ for
$S_\Q(n,d)$, $V$ for $V_\Q$.

In this paper, we give an alternative construction of Schur algebras,
as follows. First we obtain a presentation of $S(n,d)$ by generators
and relations (Theorem \ref{thm:present:S}).  This presentation is
compatible with Serre's presentation of the universal enveloping
algebra $U = U(\gl_n)$ of the Lie algebra $\gl_n$ of $n \times n$
matrices. Then we construct $S_\Z(n,d)$ as the precise analogue of the
Kostant $\Z$-form $U_\Z$, and we obtain new bases for $S_\Z(n,d)$
(Theorem \ref{thm:idemp:basis}).  We also obtain a second presentation
of $S(n,d)$ by generators and relations (Theorem
\ref{thm:idemp:present}) which is closely related to Lusztig's
construction (see \cite{Lusztig:book}) of the modified form $\dot{\U}$
of the quantized enveloping algebra $\U$.

Our approach is based on the classic ``double-centralizer'' theory of
Schur \cite{Schur} (and its quantization). The group $\GL(V)$ acts on
$V^{\otimes d}$ by means of the natural action in each tensor factor,
and by differentiating this action on tensors we obtain an action of
$U$ on $V^{\otimes d}$.  These actions obviously commute with the
action of the symmetric group $\Sigma_d$. So we have representations
$$
\begin{CD}
U @>>> \End(V^{\otimes d}) @<<< \Q\Sigma_d
\end{CD}
$$
induced from the commuting actions.  Then Schur's result is that the
image of each representation is precisely the commuting algebra for
the action of the other algebra. In particular, $S(n,d)$ is the image
of the representation $U \to \End(V^{\otimes d})$.  It is a very
natural problem to ask for an efficient generating set for the kernel
of this representation.  By solving this problem we obtain the
presentation of Theorem \ref{thm:present:S}. We note that, in the
quantum case, the analogue of the surjective map $U \to S(n,d)$ was
studied by R.M.\ Green \cite{RG:thesis, RGreen}, who described a basis
for the kernel.

In the quantum case one replaces in the above setup $\Q$ by $\Q(v)$
($v$ an indeterminate), $V$ by an $n$-dimensional $\Q(v)$-vector space
$\V$, $\Sigma_d$ by the corresponding Hecke algebra
$\mathbf{H}=\mathbf{H}(\Sigma_d)$, and $U$ by the Drinfeld-Jimbo
quantized enveloping algebra $\U = \U(\gl_n)$.  Then the resulting
commuting algebra, $\S(n,d)$, is known as the $q$-Schur algebra, or
quantized Schur algebra. It appeared first in work of Dipper and James
\cite{DJ1,DJ2}, and, independently, Jimbo \cite{Jimbo}.  Dipper and
James showed that the $q$-Schur algebras determine the representation
theory of the finite general linear groups in non-describing
characteristic. (Note that one should replace their parameter $q$ by
$v^2$ to make the correspondence with our version of $\S(n,d)$.) In
\cite{BLM} a geometric realization of $\S(n,d)$ was given.  In
\cite{Du} the \cite{BLM} approach was reconciled with the Dipper-James
approach; moreover, it was shown in that paper that $\S(n,d)$ may be
identified with the image of the map
$$
\begin{CD}
\U @>>> \End(\V^{\otimes d}).
\end{CD}
$$
In the quantum situation, one replaces $\Z$ by
$\A=\Z[v,v^{-1}]$. In this case the analogue of
$S_\Z(n,d)$ is a certain $\A$-form $\S_\A(n,d)$ in $\S(n,d)$.

Our results in the quantum case are almost exact analogues of the
results in the classical case, although the proofs are sometimes more
difficult.  First we obtain a presentation of $\S(n,d)$ by generators
and relations (Theorem \ref{qthm:present:S}).  This presentation is
compatible with the usual presentation of the quantized enveloping
algebra $\U$ (over $\Q(v)$) corresponding to $\gl_n$.  Then we
construct $\S_\A(n,d)$ as the analogue of Lusztig's $\A$-form $\U_\A$,
and we obtain new bases for $\S_\A(n,d)$ (Theorem
\ref{qthm:idemp:basis}).  Finally, we have a second presentation
of $\S(n,d)$ by generators and relations (Theorem
\ref{qthm:idemp:present}) which is closely related to the algebra
$\dot{\U}$.  Upon specializing $v$ to $1$, the presentation of Theorem
\ref{qthm:idemp:present} coincides with the presentation of Theorem
\ref{thm:idemp:present}.  (This does not apply to Theorem
\ref{qthm:present:S} in relation to Theorem \ref{thm:present:S}.)

In the final section we give some applications of our results to the
Borel Schur algebras and Hecke algebras.  In particular, in
Proposition \ref{Borelbasis} we obtain a simple basis for the Borel
Schur algebras which is a subset of the integral basis obtained (in
Theorem \ref{qthm:idemp:basis}) for the entire Schur algebra. We also
get a new basis for the Hecke algebra (realized as a subalgebra of the
Schur algebra) which is a subset of our basis in Theorem
\ref{qthm:idemp:basis}, and an integral presentation of it.
We also write out some examples in the final section.

In rank $1$ we have more precise results than in this paper, obtained
by different arguments \cite{DG}.  The results in the current
paper were summarized in the announcement \cite{DG:announce}.

{\em Acknowledgements.} The authors have benefitted from discussions with
T.\ Br\"ustle, K.\ Erdmann, R.M.\ Green, S.\ K\"onig, R.\ Marsh. 
We are grateful to the referee for useful suggestions.

\section{Main results: classical case} \label{sec:main}

Let $\Phi$ be the root system of type $A_{n-1}$: $\Phi =
\{\varepsilon_i - \varepsilon_j \mid 1 \le i \ne j \le n \}$.  Here
the $\varepsilon_i$ form the standard orthonormal basis of the euclidean
space $\R^n$.  Let $(\ ,\ )$ denote the inner product on this space
and define $\alpha_i = \varepsilon_i - \varepsilon_{i+1}$.  Then
$\{\alpha_1, \dots, \alpha_{n-1} \}$ is a base of simple roots and
$\Phi^+ = \{ \varepsilon_i - \varepsilon_j \mid i<j \}$ is the
corresponding set of positive roots.

We now give a precise statement of our main results in the classical
case. The proofs are contained in sections
\ref{sec:B}--\ref{sec:idemp:present}. The first result describes a
presentation by generators and relations of the Schur algebra over the
rational field $\Q$.

\begin{thm} \label{thm:present:S}
Over $\Q$, the Schur algebra $S(n,d)$ is isomorphic with the
associative algebra (with 1) on the generators 
$e_i$, $f_i$ ($1 \le i \le n-1$), and $H_i$ ($1 \le i \le n$) with
relations
\begin{align}
H_i H_j &= H_j H_i \tag{R1}\label{R1} \\
e_i f_j - f_j e_i &= \delta_{ij}(H_j - H_{j+1})  \tag{R2}\label{R2}
\end{align}
\begin{equation}\label{R3}
H_i e_j - e_j H_i = (\varepsilon_i, \alpha_j) e_j, \quad
H_i f_j - f_j H_i = -(\varepsilon_i, \alpha_j) f_j \tag{R3}
\end{equation}
\begin{equation}\label{R4}
\begin{aligned}
{}&e_i^2 e_j - 2 e_i e_j e_i + e_j e_i^2 = 0 \quad (|i-j|=1)\\
&e_i e_j - e_j e_i = 0 \quad (\text{otherwise})
\end{aligned} \tag{R4}
\end{equation}
\begin{equation}\label{R5}
\begin{aligned}
{}&f_i^2 f_j - 2 f_i f_j f_i + f_j f_i^2 = 0 \quad (|i-j|=1)\\
&f_i f_j - f_j f_i = 0 \quad (\text{otherwise})
\end{aligned} \tag{R5}
\end{equation}
\begin{align}
 H_1 + H_2 + \cdots + H_n & = d \tag{R6}\label{R6} \\
 H_i(H_i-1)\cdots(H_i-d) & = 0. \tag{R7}\label{R7}
\end{align}
\end{thm}
Note that the enveloping algebra $U=U(\gl_n)$ is the algebra on the
same generators but subject only to the relations
\eqref{R1}--\eqref{R5}, and $U(\sl_n)$ is isomorphic with the
subalgebra of $U$ generated by the $e_i, f_i, H_i - H_{i+1}$ ($1 \le i
\le n-1$).

Next we introduce the root vectors $x_\alpha$ ($\alpha\in
\Phi$), which may be defined inductively as follows. Write
$\alpha=\varepsilon_i-\varepsilon_j$ and assume that $i<j$. If $j-i =
1$ then $\alpha=\alpha_i$ and we set $x_\alpha=e_i$,
$x_{-\alpha}=f_i$.  If $j-i>1$ then we inductively set
$$
 x_\alpha = e_ix_{\alpha-\alpha_i} - x_{\alpha-\alpha_i}e_i, \qquad 
 x_{-\alpha} = x_{-\alpha+\alpha_i} f_i - f_i x_{-\alpha+\alpha_i}.
$$
The linear span of the set $\{ x_\alpha \} \cup \{ H_i \}$ is a
subspace of $U$ isomorphic with the Lie algebra $\gl_n$ under the Lie
bracket given by $[x,y]=xy-yx$, and the $x_\alpha$ correspond to the
usual root vectors in $\gl_n$.

\begin{remark} \label{thm:truncPBW}
The defining relation \eqref{R6} can be used to rewrite one of
the $H_i$'s in terms of the others. 
Fix an integer $i_0$ with $1\le i_0 \le n$ and set
$$
G=\{ x_\alpha \mid \alpha\in \Phi \} \cup
\{H_i \mid i \ne i_0 \}
$$
and fix an arbitrary ordering for this set.  We conjecture that
$S(n,d)$ has a $\Q$-basis consisting of all monomials in $G$
(with specified order) of total degree not exceeding $d$. This basis
would be an analogue of the Poincare-Birkhoff-Witt (PBW) basis of $U$.
\end{remark}

Our next result constructs the integral Schur algebra $S_\Z(n,d)$
in terms of the generators given above.  We need more
notation. For $B=(B_i)$ in $\N^n$, we write
$$
H_B = \prod_{i=1}^n \dbinom{H_i}{B_i}
$$
where $\dbinom{H_i}{m} = H_i(H_i-1)\dots (H_i-m+1)/(m!)$ ($m \ge 1$),
$\dbinom{H_i}{0} = 1$.  Let $\Lambda(n,d)$ be the subset of $\N^n$
consisting of those $\lambda \in \N^n$ satisfying $|\lambda|=d$ (here
$|\lambda|{}=\sum \lambda_i$); this is the set of $n$-part
compositions of $d$. Given $\lambda \in \Lambda(n,d)$ we set
$1_\lambda{}=H_{\lambda}.$ We will show that the collection
$\{1_\lambda\}$ as $\lambda$ varies over $\Lambda(n,d)$ forms a set of
pairwise orthogonal idempotents in $S_\Z(n,d)$ which sum to the
identity element.

For $m \in \N$ and $\alpha \in \Phi$, set $\divided{x_\alpha}{m}{}=
x_\alpha^m/(m!)$. Any product of elements of the form
$$
\divided{x_\alpha}{r}, \quad \dbinom{H_i}{s}\qquad
(r,s \in \N, \alpha\in\Phi, 1\le i \le n),
$$
taken in any order, will be called a {\em Kostant monomial}.  Note
that the set of Kostant monomials is multiplicatively closed.  We
define a function $\ct$ (content function) on Kostant monomials by
setting
$$
\ct(\divided{x_\alpha}{m}){}= m\,\varepsilon_{\max(i,j)}, \quad
\ct( \binom{H_i}{m} ) {}= 0
$$
where $\alpha = \varepsilon_i - \varepsilon_j$ ($i\ne j$), and by
declaring that $\ct(XY) = \ct(X)+\ct(Y)$ whenever $X,Y$ are Kostant
monomials.

For $A\in \N^{\Phi^+}$ 
we set $|A|{}= \sum_{\alpha\in
\Phi^+} A(\alpha)$.  For $A, C \in \N^{\Phi^+}$ we write
$$
e_A = \prod_{\alpha\in \Phi^+} \divided{x_\alpha}{A(\alpha)}, \quad
f_C = \prod_{\alpha\in \Phi^+} \divided{x_{\alpha}}{C(-\alpha)}
$$
where the products in $e_A$ and $f_C$ are taken relative to any
two fixed orders on $\Phi^+$.

The first part of the next result shows that $S_\Z(n,d)$ is the
analogue in $S(n,d)$ of Kostant's $\Z$-form $U_\Z$ in $U$.

\begin{thm} \label{thm:idemp:basis} The integral Schur
algebra $S_\Z(n,d)$ is the subring of $S(n,d)$ generated by all
divided powers $\divided{e_i}{m},\, \divided{f_i}{m}$.  Moreover, each
of the disjoint unions
\begin{align}
\basis_+ &= \textstyle\bigcup_{\lambda}\, 
 \{e_A 1_\lambda f_C \mid \ct(e_A f_C) \preceq \lambda\}\tag{a}\\
\basis_-&= \textstyle\bigcup_{\lambda}\, 
 \{f_A 1_\lambda e_C \mid \ct(f_A e_C) \preceq \lambda\},\tag{b}
\end{align}
as $\lambda$ varies over $\Lambda(n,d)$, and
where $\preceq$ denotes the {\em componentwise} partial ordering on
$\N^n$, is a $\Z$-basis of $S_\Z(n,d)$.
\end{thm}

Finally, we have another presentation of the Schur algebra by
generators and relations. This presentation has the advantage that it
possesses a quantization of the same form, in which we can specialize
$v$ to $1$ to recover the classical version.

\begin{thm} \label{thm:idemp:present}
The $\Q$-algebra $S(n,d)$ is the associative
algebra (with $1$) given by generators $1_\lambda$ ($\lambda\in
\Lambda(n,d)$), $e_i$, $f_i$ ($1\le i \le n-1$) subject to
the relations
\begin{equation} \label{S1}
1_\lambda 1_\mu = \delta_{\lambda,\mu} 1_\lambda, \quad
\sum_{\lambda\in \Lambda(n,d)} 1_\lambda = 1 \tag{R$1^\prime$}
\end{equation}
\begin{equation}\label{S2}
\begin{gathered}
e_i 1_\lambda =
\begin{cases}
1_{\lambda+\alpha_i} e_i &
   \text{if $\lambda+\alpha_i \in \Lambda(n,d)$}\\
0 & \text{otherwise}
\end{cases}\\
f_i 1_\lambda =
\begin{cases}
1_{\lambda-\alpha_i} f_i &
   \text{if $\lambda-\alpha_i \in \Lambda(n,d)$}\\
0 & \text{otherwise}
\end{cases}\\
1_\lambda e_i =
\begin{cases}
e_i 1_{\lambda-\alpha_i} &
   \text{if $\lambda-\alpha_i \in \Lambda(n,d)$}\\
0 & \text{otherwise}
\end{cases}\\
1_\lambda f_i =
\begin{cases}
f_i 1_{\lambda+\alpha_i} &
   \text{if $\lambda+\alpha_i \in \Lambda(n,d)$}\\
0 & \text{otherwise}
\end{cases}
\end{gathered}\tag{R$2^\prime$}
\end{equation}
\begin{equation} \label{S3}
e_i f_j - f_j e_i = \delta_{ij}
\sum_{\lambda\in \Lambda(n,d)} (\lambda_j-\lambda_{j+1}) 1_\lambda \tag{R$3^\prime$}
\end{equation}
along with the Serre relations \eqref{R4}, \eqref{R5}.
\end{thm}

\section{Main results: quantum case} \label{qsec:main}

Our main results in the quantum case are similar in form to those in
the classical case. Proofs are given in sections
\ref{qsec:B}--\ref{qsec:truncPBW}.  The first result describes a
presentation by generators and relations of the quantized Schur
algebra over the rational function field $\Q(v)$.

\begin{thm} \label{qthm:present:S}
Over $\Q(v)$, the $q$-Schur algebra $\S(n,d)$ is isomorphic with the
associative algebra (with 1) with 
generators $E_i$, $F_i$ ($1\leq i\leq n-1$), $K_i$, $K_i^{-1}$
($1\leq i\leq n$) and relations
\begin{align}
& K_i K_j =K_j K_i, \qquad K_i K_i^{-1} = K_i^{-1}K_i =1
  \tag{Q1}\label{Q1}\\
& E_iF_j - F_jE_i =
\delta_{ij}\frac{K_iK_{i+1}^{-1}-K_i^{-1}K_{i+1}}{v-v^{-1}}\tag{Q2}\label{Q2}\\
& K_i E_j = v^{(\ep_i, \al_j)}E_jK_i, \qquad
K_i F_j = v^{-(\ep_i, \al_j)}F_jK_i \tag{Q3}\label{Q3}
\end{align}
\begin{equation}\label{Q4}
\begin{aligned}
{}& E_i^2E_j - (v+v^{-1})E_iE_jE_i +E_jE_i^2 =0 \quad(|i-j|=1) \\
{}& E_iE_j - E_i E_j = 0 \quad(\text{otherwise})
\end{aligned}\tag{Q4}
\end{equation}
\begin{equation}
\begin{aligned}\label{Q5}
{}& F_i^2F_j - (v+v^{-1})F_iF_jF_i +F_jF_i^2 =0 \quad(|i-j|=1) \\
{}& F_iF_j - F_jF_i = 0 \quad(\text{otherwise})
\end{aligned}\tag{Q5}
\end{equation}
\begin{gather}
 K_1  K_2 \cdots  K_n = v^d \tag{Q6}\label{Q6} \\
 (K_i-1)(K_i-v)(K_i-v^2)\cdots(K_i-v^d) = 0. \tag{Q7}\label{Q7}
\end{gather}
\end{thm}

We note that the quantized enveloping algebra $\U=\U(\gl_n)$ is the
algebra on the same set of generators, subject only to the relations
\eqref{Q1}--\eqref{Q5}. Moreover, the quantized enveloping algebra
$\U(\sl_n)$ is isomorphic with the subalgebra of $\U$ generated by all
$E_i, F_i, K_iK_{i+1}^{-1}$ ($1\le i \le n-1$).

We have $q$-analogues of the root vectors in $\gl_n$, which can be
defined as follows.  For $\alpha\in \Phi^+$, write
$\alpha=\varepsilon_i-\varepsilon_j$ for $i < j$.  If $j-i=1$ then set
$X_\alpha{}=E_i$, $X_{-\alpha}{}=F_i$. For $j-i>1$ we inductively set
(following Xi \cite[\sect5.6]{Xi})
$$ 
X_\alpha{}=v^{-1}E_i X_{\alpha-\alpha_i} - X_{\alpha-\alpha_i}E_i,
\quad X_{-\alpha}{}=vX_{-\alpha+\alpha_i}F_i - F_i
X_{-\alpha+\alpha_i}.
$$
Our notation differs from that in \cite{Xi} where the
elements $X_\alpha$ and $X_{-\alpha}$ are denoted $E_{i,j-1}$
and $F_{j,i+1}$, resp.  Up to scalar multiplication by units in $\A$,
the elements $X_\alpha$ ($\alpha\in \Phi$) first appeared in
Jimbo's paper \cite{Jimbo}.

\begin{remark} \label{qthm:truncPBW} 
Fix an integer $i_0$ with $1\le i_0 \le n$ and write $\N^n_{i_0}$ for
the set of $B \in \N^n$ such that $B_{i_0} = 0$.  We conjecture that
$\S(n,d)$ has a  $\Q(v)$-basis consisting of all monomials of the form
$$
\prod_{\alpha\in \Phi^+} X_\alpha^{A(\alpha)} \prod_{i\ne i_0} K_i^{B_i} 
\prod_{\alpha\in \Phi^+} X_{-\alpha}^{C(\alpha)}
 \qquad (A,C \in \N^{\Phi^+}, B \in \N^n_{i_0})
$$
of total degree not exceeding $d$, where the products of powers of
$X_\alpha$, $X_{-\alpha}$ are taken with respect to arbitrary fixed
orders on $\Phi^+$.  This basis would be the analogue of the PBW-type
basis of $\U$, given in Lusztig \cite[Proposition 1.13]{Lusztig}.
\end{remark}

Our next result constructs the $\A$-form $\S_\A(n,d)$ in terms of the
generators given above. We write $\sqbinom{K_i}{t}$ short for
$\sqbinom{K_i;0}{t}$, where (following Lusztig) we define
$$
\sqbinom{K_i; c}{t}{}=\prod_{s=1}^{t}
\frac{K_i v^{c-s+1} - K_i^{-1} v^{-c+s-1}}{v^s-v^{-s}}, 
\qquad \sqbinom{K_i; c}{0} = 1
$$  
for $t\ge 1$, $c\in \Z$.  For $B$ in $\N^n$, we write
$$
K_B = \prod_{i=1}^n \sqbinom{K_i}{B_i}.
$$
Given $\lambda \in \Lambda(n,d)$ we set $1_\lambda{}= K_\lambda$. 
Just as in the classical case, the
collection $\{1_\lambda\}$ forms a set of pairwise orthogonal
idempotents in $\S_\A(n,d)$ which sum to the identity element.

For $m \in \Z$ let $[m]$ denote the quantum integer $[m]{}=
(v^m-v^{-m})/(v-v^{-1})$ and set
\begin{gather*}
[m]! {}= [m][m-1] \cdots [1], \qquad [0]! = 1 \\
\sqbinom{c}{m} = \frac{[c][c-1]\cdots[c-m+1]}{[m]!}, \qquad 
\sqbinom{c}{0} = 1
\end{gather*}
for $c\in \Z$, $m \ge 1$.  The $q$-analogues of the divided powers of
root vectors are defined by $\divided{X_\alpha}{m}{}= X_\alpha/[m]!$.
The Kostant monomials in this situation are products of elements of
the form
$$
\divided{X_\alpha}{r}, \quad\sqbinom{K_i}{s},
\quad K_i^{\pm 1}
\qquad (r, s \in \N, \alpha\in\Phi, 1\le i \le n),
$$
taken in any order.  As before, the set of Kostant monomials is
multiplicatively closed.  By analogy with the classical case, $\ct$ is
defined by
$$
\ct(\divided{X_\alpha}{m})
  {}= m\,\varepsilon_{\max(i,j)}, \quad
\ct( \sqbinom{K_i}{m} ) = \ct( K_i^{\pm 1} ) {}= 0
$$
where $\alpha = \varepsilon_i - \varepsilon_j \in \Phi$, and by
declaring that $\ct(XY) = \ct(X)+\ct(Y)$ whenever $X,Y$ are
Kostant monomials. For $A, C \in \N^{\Phi^+}$ we write
$$
E_A = \prod_{\alpha\in \Phi^+} \divided{X_\alpha}{A(\alpha)}, \quad
F_C = \prod_{\alpha\in \Phi^+} \divided{X_{-\alpha}}{C(\alpha)}
$$
where the products in $E_A$ and $F_C$ are taken relative to any
two specified orderings on $\Phi^+$.

\begin{thm} \label{qthm:idemp:basis}
The integral $q$-Schur algebra $\S_\A(n,d)$ is the subring of
$\S(n,d)$ generated by all quantum divided powers 
$\divided{E_i}{m}, \divided{F_i}{m}$, 
along with the elements $\sqbinom{K_i}{m}$.
Moreover, each of the sets
\begin{align}
\qbasis_+ &= \textstyle\bigcup_\lambda\, 
 \{E_A 1_\lambda F_C \mid \ct(E_A F_C) \preceq \lambda\}\tag{a}\\
\qbasis_- &= \textstyle\bigcup_\lambda\, 
 \{F_A 1_\lambda E_C \mid \ct(F_A E_C) \preceq \lambda\},\tag{b}
\end{align}
as $\lambda$ ranges over $\Lambda(n,d)$, forms an $\A$-basis of
$\S_\A(n,d)$.
\end{thm}

We conjecture that the elements $\sqbinom{K_i}{m}$ lie within the subring
generated by the $\divided{E_i}{m}$, $\divided{F_i}{m}$, in which case
we would obtain the more precise analogue of Theorem \ref{thm:idemp:basis}.

Finally, we have another presentation of the $q$-Schur algebra by
generators and relations. These relations are similar to relations
that hold for the modified form $\dot{\U}$ of $\U$ (see \cite[Chap.\
23]{Lusztig:book}).  This presentation has the advantage that upon
specializing $v$ to $1$, we recover the classical version given in
Theorem \ref{thm:idemp:present}.

\begin{thm} \label{qthm:idemp:present}
The algebra $\S(n,d)$ is the associative
algebra (with $1$) given by generators $1_\lambda$ ($\lambda\in
\Lambda(n,d)$), $E_i$, $F_i$ ($1 \le i \le n-1$) subject to
the relations
\begin{equation}\label{S1'}
1_\lambda 1_\mu = \delta_{\lambda,\mu} 1_\lambda, \quad
\sum_{\lambda\in \Lambda(n,d)} 1_\lambda = 1 \tag{Q$1^\prime$}
\end{equation}
\begin{equation}\label{S2'}
\begin{gathered}
E_i 1_\lambda =
\begin{cases}
1_{\lambda+\alpha_i} E_i &
   \text{if $\lambda+\alpha_i \in \Lambda(n,d)$}\\
0 & \text{otherwise}
\end{cases}\\
F_i 1_\lambda =
\begin{cases}
1_{\lambda-\alpha_i} F_i &
   \text{if $\lambda-\alpha_i \in \Lambda(n,d)$}\\
0 & \text{otherwise}
\end{cases}\\
1_\lambda E_i =
\begin{cases}
E_i 1_{\lambda-\alpha_i} &
   \text{if $\lambda-\alpha_i \in \Lambda(n,d)$}\\
0 & \text{otherwise}
\end{cases}\\
1_\lambda F_i =
\begin{cases}
F_i 1_{\lambda+\alpha_i} &
   \text{if $\lambda+\alpha_i \in \Lambda(n,d)$}\\
0 & \text{otherwise}
\end{cases}
\end{gathered}\tag{Q$2^\prime$}
\end{equation}
\begin{equation} \label{S3'}
E_i F_j - F_j E_i = \delta_{ij}
\sum_{\lambda\in \Lambda(n,d)} [\lambda_j-\lambda_{j+1}] 1_\lambda
  \tag{Q$3^\prime$}
\end{equation}
along with the $q$-Serre relations \eqref{Q4}, \eqref{Q5}.
\end{thm}

\section{The algebra $\B$}\label{sec:B}

From now on we hold $n$ and $d$ fixed, and set $S{}=S(n,d)$.  We define
an algebra $\B=\B(n,d)$ (over $\Q$) by the generators and relations of
Theorem \ref{thm:present:S}.  Since $U$ is the algebra on the same
generators but subject only to relations \eqref{R1}--\eqref{R5}, we
have a surjective quotient map $U \to \B$ (mapping generators onto
generators).  Eventually we shall show that $\B \simeq S$, which will
prove Theorem \ref{thm:present:S}.

\begin{lem}\label{lem:minpoly}
Under the representation $U \to \End(V^{\otimes d})$ the images of the
$H_i$ satisfy the relations \eqref{R6} and \eqref{R7}.
Moreover, the relation \eqref{R7} is the minimal polynomial
of (the image of) $H_i$ in $\End(V^{\otimes d})$.
\end{lem}

\begin{proof}
Relation \eqref{R6} is trivial in the case $d=1$, from which the
general case follows since each $H_i$ acts as a derivation of
$V^{\otimes d}$. The relation \eqref{R7} follows from the fact (which
can be verified by induction on $d$) that the eigenvalues of the
diagonal operators $H_i$ are $0,1,\ldots,d$.  The proof is complete.
\end{proof}

As we know, $S = S(n,d)$ is the image of the representation $U \to
\End(V^{\otimes d})$ mentioned in the introduction.  From the above
lemma it follows that this surjection $U \to S$ factors through
$\B$.  Because $\B$, $S$ are homomorphic images of $U$, any relations
between generators holding in $U$ will automatically carry over to
$\B$, $S$. We will not distinguish notationally between the generators
or root vectors for $U$, $\B$, or $S$.

Recall the triangular decomposition of $U$, that the multiplication map
$U^- \otimes U^0 \otimes U^+ \xrightarrow{\approx} U$
is an isomorphism of vector spaces, where $U^+$ (resp., $U^-$) is the
subalgebra of $U$ generated by the $e_i$ (resp., $f_i$) and $U^0$ is
the subalgebra of $U$ generated by all $H_i$. 
Thus $U = U^- U^0 U^+$.  From
this we obtain a similar triangular decomposition of $\B$:
\begin{equation}\label{B:TD}
\B = \B^- \B^0 \B^+
\end{equation}
where $\B^+, \B^-, \B^0$ are defined to be the images of $U^+,
U^-, U^0$ under the quotient mapping $U \to \B$.

We also have similar factorizations over $\Z$.  Setting $U_\Z^+$,
$U_\Z^-$, $U_\Z^0$ to be, respectively, the intersection of $U^+$,
$U^-$, $U^0$ with the Kostant $\Z$-form $U_\Z$ (the subring of $U$
generated by all $\divided{e_i}{m}$, $\divided{f_i}{m}$,
$\dbinom{H_i}{m}$), we have the factorization $ U_\Z = U_\Z^- U_\Z^0
U_\Z^+ $, which immediately induces similar equalities
\begin{equation}\label{UZ:TD}
\B_\Z = \B_\Z^- \B_\Z^0 \B_\Z^+
\end{equation}
where the various subalgebras are defined in the obvious way as
appropriate homomorphic images of $U_\Z^+$, $U_\Z^-$, $U_\Z^0$.

Since $U_\Z^+$ (resp., $U_\Z^-$) is the $\Z$-subalgebra of $U$
generated by the $\divided{x_\alpha}{m}$ for $\alpha\in \Phi^+$
(resp., $\alpha \in \Phi^-$) and $m\in \N$, the same statement applies
to $\B_\Z^+$ (resp., $\B_\Z^-$) in relation to $\B$.  Moreover,
$U_\Z^0$ is the $\Z$-subalgebra of $U$ generated by the
$\dbinom{H_i}{m}$ for $1\le i \le n$ and $m\in \N$, so $\B_\Z^0$ is
the $\Z$-subalgebra of $\B$ generated by the same elements.

Now we investigate the structure of the algebra $\B^0$.  We start with
the algebra $U^0$, which is isomorphic with the polynomial ring
$\Q[H_1,\dots,H_n]$ in $n$ commuting indeterminates $H_1,\dots, H_n$.
By the remarks following Lemma \ref{lem:minpoly} we have surjections
$U \to \B \to S$.  Let $S^0$ be the image in $S$ of $U^0$ under the
map $U \to S$. Clearly we have surjections $U^0 \to \B^0 \to S^0$
obtained from $U \to \B \to S$ by restriction.

\begin{prop} \label{prop:idempotent}
Define an algebra $\B^\prime = U^0/I^0$ where $I^0$ is the ideal in
$U^0$ generated by elements $H_i(H_i-1)\cdots(H_i-d)\ \ (1\le i \le
n)$ and $H_1+\cdots+H_n-d$.

\par\noindent(a) We have an algebra isomorphism $\B^\prime \cong \B^0$. 

\par\noindent(b) The set $\{ 1_\lambda \mid \lambda \in \Lambda(n,d)
\}$ is a $\Q$-basis for $\B^0$ and a $\Z$-basis for $\B_\Z^0$;
moreover, this set is a set of pairwise orthogonal idempotents which
add up to $1$.

\par\noindent(c) $H_B = 0$ for any $B\in \N^n$
such that $|B|>d$.
\end{prop}

\begin{proof}
\newcommand{\tB}{\widetilde{\B}}
Consider first the algebra $\tB^\prime$ defined to be the quotient of
$U^0$ by the ideal generated only by the elements
$H_i(H_i-1)\cdots(H_i-d)$ $(1 \le i \le n)$. 
Since each of the relations in $\tB^\prime$ is a polynomial in just one 
of the variables, we have the factorization
$$
\tB^\prime\cong \Q [H_1]/(p(H_1)) \otimes \cdots \otimes \Q [H_n]/(p(H_n))
$$
where $p(X)=X(X-1)\cdots(X-d)$.  By the Chinese Remainder Theorem
applied to each tensor factor we obtain from the above isomorphisms
(products denote direct products)
\begin{align*}\tB^\prime&\cong
\prod_{i=0}^d \left( \Q [H_1]/(H_1 -i)\right) \otimes \cdots
\otimes \prod_{i=0}^d \left( \Q [H_n]/(H_n -i)\right)\\ 
\intertext{and by rearranging the order of factors we obtain}
& \cong \prod_{0\leq \mu_1,\ldots, \mu_n \leq d}\left( \Q [H_1]/(H_1-\mu_1)
\otimes \cdots \otimes \Q [H_n]/(H_n-\mu_n)\right)\\ 
& \cong \prod_{0\leq \mu_1,\ldots, \mu_n \leq d} \Q [H_1, \ldots,
H_n]/(H_1-\mu_1, \ldots, H_n - \mu_n).
\end{align*}
The isomorphism is realized by the map which sends a polynomial
$f(H_1,\ldots, H_n)$ in $\tB^\prime$ to the element $( f(\mu_1, \ldots, \mu_n)
)_{0\leq \mu_1,\ldots, \mu_n \leq d}$ of the direct product.
Since $\B^\prime$ is isomorphic with $\tB^\prime/(H_1+\cdots+H_n-d)$, we
deduce from the above that
$$
\B^\prime \cong \prod_{\mu \in \Lambda(n,d)}\Q [H_1, \ldots,
H_n]/(H_1-\mu_1, \ldots, H_n - \mu_n)
$$
and this isomorphism (which we denote by $\phi$) is realized by the
map sending $f(H_1,\ldots, H_n)$ to the element $( f(\mu_1, \ldots,
\mu_n) )_{\mu \in \Lambda(n,d)}$.

Thus, given any $\lambda \in \Lambda(n,d)$, we have
$$
\phi(1_\lambda) = \left(\, \dbinom{\mu_1}{\lambda_1} \cdots
\dbinom{\mu_n}{\lambda_n}\, \right)_{\mu \in \Lambda(n,d)} =
( \delta_{\lambda\mu} )_{\mu \in \Lambda(n,d)}
$$
Since the vector on the right of the above equality consists of zeros
and ones, with precisely one nonzero entry, and since $\phi$ is an
isomorphism, it follows that the various $1_{\lambda}$ are orthogonal
idempotents which add up to the identity of $\B^\prime$. 

Now let $I$ be the ideal in $U$ generated by elements
$H_i(H_i-1)\cdots(H_i-d)\ \ (1\le i \le n)$ and $H_1+\cdots+H_n-d$. 
Then by definition $\B \cong U/I$.  The canonical quotient map
$U \to U/I$ induces, upon restriction to $U^0$, a map $U^0 \to U/I$.
The image of this map is $\B^0$ and its kernel is $U^0 \cap I$, so
$\B^0 \cong U^0/(U^0 \cap I)$.  Clearly $I^0 \subset U^0 \cap I$. 
Thus we obtain the following sequence of algebra surjections
\begin{equation}
\begin{CD}
\B^\prime = U^0/I^0 @>>> U^0/(U^0 \cap I) @>\approx>> \B^0 @>>> S^0
\end{CD}
\end{equation}
where the last map is obtained by Lemma \ref{lem:minpoly}, and the
middle map is actually an isomorphism. By the above we see that the
dimension of $\B^\prime$ is the cardinality of the set $\Lambda(n,d)$. This
is well-known to equal the dimension of the zero part $S^0$ of the Schur
algebra $S=S(n,d)$. It follows that all the surjections above are
algebra isomorphisms.  This proves parts (a) and (b).

To prove part (c), suppose that $B\in \N^n$ such that $|B|>d$. Then
for each $\mu \in \Lambda (n,d)$, there exists $i$ with
$\mu_i<b_i$. Thus $\phi(H_B)=0$ and $H_B$ must be $0$ since the map
$\phi$ is an isomorphism. The proof is complete.
\end{proof}

The next result is obtained by a similar argument. 

\begin{prop} \label{prop:Hmult}
Let $1\le i \le n$, $b\in \N$, $\lambda \in \Lambda(n,d)$, and
$B \in \N^n$.  We have the following identities in
the algebra $\B^0$:
\begin{gather} H_i 1_\lambda = \lambda_i 1_\lambda, \qquad
\dbinom{H_i}{b} 1_\lambda = \dbinom{\lambda_i}{b} 1_\lambda\tag{a}\\
H_B \, 1_\lambda = \lambda_B\,1_\lambda, \quad
\text{where $\lambda_B {}= \prod_i \dbinom{\lambda_i}{B_i}$}\tag{b}\\
H_B = \sum_{\lambda} \lambda_B \, 1_\lambda\tag{c}
\end{gather}
where the sum in part (c) is carried out over all $\lambda \in
\Lambda(n,d)$.
\end{prop}

\begin{proof} We apply the isomorphism $\phi$ from the proof of the
 preceding proposition to the product on the left-hand-side of (a),
obtaining
$$
\phi(\binom{H_i}{b}\, 1_\lambda) = \phi(\binom{H_i}{b}
\binom{H_1}{\lambda_1} \cdots \binom{H_n}{\lambda_n}).
$$
The right-hand side of the above equality gives the vector
$$
\left( \, \binom{\mu_i}{b} \binom{\mu_1}{\lambda_1} \cdots
\binom{\mu_n}{\lambda_n} \, \right)_{\mu\in\Lambda(n,d)}
= \left( \, \binom{\mu_i}{b} \delta_{\mu\lambda} \, \right)_{\mu
\in\Lambda(n,d)} 
$$
which is the same as $\dbinom{\lambda_i}{b} \phi(1_\lambda)$ since
from the preceding proof $\phi(1_\lambda) = (\delta_{\mu\lambda})_\mu$.
Since $\phi$ is an isomorphism, part (a) is proved.

Part (b) follows immediately from part (a).  Then by the result of
part (b) we obtain the equalities $H_B = H_B \cdot 1
= H_B \sum_\lambda 1_\lambda = \sum_\lambda \lambda_B 1_\lambda$,
proving part (c).
\end{proof}

We write $\N^n_{i_0}$ for the set of $B=(B_i)\in \N^n$ such that
$B_{i_0}=0$.

\begin{cor}\label{cor:H_B-basis}
For any fixed choice of $i_0$ ($1\le i_0 \le n$) the
set $\{ H_B \mid B\in \N^n_{i_0}, |B| \le d \}$ is a $\Q$-basis for
$\B^0$ and a $\Z$-basis for $\B_\Z^0$.
\end{cor}

\begin{proof}

The sets $\Lambda(n,d)$ and $\{B\in \N_{i_0}^n,\,\, |B|\leq d\}$ have
the same cardinality. For instance, the map
$$
  \lambda \to \lambda - \lambda_{i_0}\varepsilon_{i_0}
$$
is bijective between the sets in question, with inverse map 
$$
   B \to B + (d-|B|)\varepsilon_{i_0}.
$$
Thus, to prove the result it is enough to show that
the set $\{H_B \mid B\in \N_{i_0}^n,\,\, |B|\leq d\}$ spans
$T_{\Z}^0$. This can be deduced by considering a certain order
(depending on $i_0$) on each of these sets.

Fixing $1\leq i_0\leq n$, we order the set $\Lambda(n,d)$ by declaring
that $\lambda$ precedes $\lambda'$ if $\lambda_{i_0}< \lambda'_{i_0}$,
or if $\lambda_{i_0}= \lambda'_{i_0}$ and there exists an index $l\ne
i_0$ such that $\lambda_j \geq \lambda_j'$ for all $j\in \{1,\ldots
l\}-\{i_0\}$. Similarly, we order the set $\{B\in \N_{i_0}^n,\,\,
|B|\leq d\}$ by declaring that $B$ precedes $B'$ if there exists an
index $l\ne i_0$ such that $b_j \geq b_j'$ for all $j\in \{1,\ldots
l\}-\{i_0\}$. With these orderings, it follows from part (c) of the
preceding proposition that the matrix of coefficients obtained by
expressing the $H_B$'s in terms of the $1_\lambda$'s is lower
unitriangular. To see this, observe that for any given $B$, with
corresponding $\lambda = B + (d-|B|)\varepsilon_{i_0}$, any $\mu$
which succeeds $\lambda$ in the above order satisfies $\mu_B = 0$, and
moreover $\lambda_B = 1$.  It follows that these equations can be
inverted over $\Z$ and so every $1_{\lambda}$ is expressible as a
$\Z$-linear combination of elements from the set $\{H_B \mid B\in
\N_{i_0}^n,\,\, |B|\leq d\}$. This proves that this set spans
$T_{\Z}^0$ and, by our remarks above, it must therefore be a
$\Q$-basis, and thus is linearly independent over $\Z$, and hence also
a $\Z$-basis. The proof is complete.
\end{proof}

The next step is to find  spanning sets for the plus part $\B^+$ and
minus part $\B^-$ of $\B$. For this we use the following result.

\begin{prop} \label{prop:rootXidemp}
For any $\alpha \in \Phi$, $\lambda \in \Lambda(n,d)$
we have the commutation formulas
$$
x_\alpha 1_\lambda =
\begin{cases}
1_{\lambda+\alpha} x_\alpha 
& \text{if $\lambda+\alpha \in \Lambda(n,d)$}\\
0 & \text{otherwise}
\end{cases}
$$
and similarly
$$
1_\lambda x_\alpha =
\begin{cases}
x_\alpha 1_{\lambda-\alpha} 
& \text{if $\lambda-\alpha \in \Lambda(n,d)$}\\
0 & \text{otherwise}.
\end{cases}
$$
\end{prop}

\begin{proof}
Write $\alpha =\varepsilon_i -\varepsilon_j$ with $i\ne j$.  
From the defining relation \eqref{R3} and the definition 
(see \sect\ref{sec:main}) of $x_\alpha$ we have
\begin{equation}\label{commute:xH}
H_l x_\alpha  = x_\alpha (H_l + (\varepsilon_l,\alpha)) .
\end{equation}
From this we obtain equalities
$$
x_{\alpha}1_{\lambda}=x_{\alpha}\dbinom{H_1}{\lambda_1}\cdots
\dbinom{H_n}{\lambda_n}=\left(\dbinom{H_i-1}{\lambda_i}
\dbinom{H_j+1}{\lambda_j}\prod_{l\neq
i,j}\dbinom{H_l}{\lambda_l}\right)x_{\alpha}.
$$
Multiplying on the left by $H_i/(\lambda_i+1)$ and then commuting with
$x_\alpha$ yields the equality
\begin{align*}
x_{\alpha} \frac{H_i+1}{\lambda_i+1}1_\lambda &= \frac{H_i}{\lambda_i+1}
\left(\dbinom{H_i-1}{\lambda_i} \dbinom{H_j+1}{\lambda_j}\prod_{l\neq
i,j}\dbinom{H_l}{\lambda_l}\right)x_{\alpha}\\
\intertext{which by Proposition \ref{prop:Hmult}(a) simplifies to give}
x_{\alpha}1_\lambda &= \left(\dbinom{H_i}{\lambda_i+1}
\dbinom{H_j+1}{\lambda_j}\prod_{l\neq
i,j}\dbinom{H_l}{\lambda_l}\right)x_{\alpha}\\
\intertext{which (if $\lambda_j > 0$) can be rewritten in the form}
x_{\alpha}1_\lambda &= \dbinom{H_i}{\lambda_i+1}\left(\dbinom{H_j}{\lambda_j} +
\dbinom{H_j}{\lambda_j-1}\right) \prod_{l\neq i,j} 
\dbinom{H_l}{\lambda_l} x_{\alpha}.
\end{align*}
The first summand on the right-hand-side of the preceding equality
vanishes, by Proposition \ref{prop:idempotent}(c). This proves the
first part of the proposition in the case $\lambda_j>0$. In case
$\lambda_j=0$ the right-hand-side is zero. This proves the first part
of the proposition.  The proof of the second part is similar.
\end{proof}

\begin{cor}\label{cor:root:nil}
We have in $\B$ the equalities $x_\alpha^{d+1} = 0$ for all $\alpha\in
\Phi$.
\end{cor}

\begin{proof}
By iterating the result of the preceding proposition we see that
$x_\alpha^{d+1} 1_\lambda = 0$ for any $\alpha\in \Phi$ and any
$\lambda \in \Lambda(n,d)$, since it is clear that $\lambda +
(d+1)\alpha$ does not belong to $\Lambda(n,d)$.  Thus we have
equalities
$$
x_\alpha^{d+1} = x_\alpha^{d+1}\cdot 1
= x_\alpha^{d+1} \sum_{\lambda \in \Lambda(n,d)} 1_\lambda
= \sum_{\lambda \in \Lambda(n,d)} x_\alpha^{d+1} 1_\lambda = 0
$$
and this proves the claim.
\end{proof}

\section{Straightening}

We need the following variants of the notion of content.  We define
functions $\ct_L$ and $\ct_R$ (left and right content) on Kostant
monomials by setting
\begin{equation}
\begin{gathered}
\ct_R(\divided{x_\alpha}{m}){}= m\,\varepsilon_j, \quad
\ct_L(\divided{x_\alpha}{m}){}= m\,\varepsilon_i, \\
\ct_R( \binom{H_i}{m} ) = \ct_L( \binom{H_i}{m} ) {}= 0
\end{gathered}
\end{equation}
where $\alpha = \varepsilon_i - \varepsilon_j$, and again using the
rule $\ct_L(XY)=\ct_L(X)+\ct_L(Y)$ (similarly for $\ct_R$) whenever
$X$ and $Y$ are Kostant monomials.  Note that for $A,C\in \N^{\Phi^+}$
we have
\begin{equation}\label{content:firstrel}
\ct(e_A) = \ct_R(e_A), \quad
\ct(f_C) = \ct_L(f_C).
\end{equation}
From this it follows immediately that
\begin{equation}\label{content:secondrel}
\ct(e_A f_C) = \ct_R(e_A) + \ct_L(f_C).
\end{equation}

\newcommand{\lambdaleft}{{\lambda^\prime}}
\newcommand{\lambdaright}{{\lambda^{\prime\prime}}}
From Proposition \ref{prop:rootXidemp} it follows that for any $A,C\in
\N^{\Phi^+}$, $\lambda \in \Lambda(n,d)$ we have equalities
\begin{equation}
e_A 1_\lambda f_C
= 1_\lambdaleft e_A f_C = e_A f_C 1_\lambdaright
\end{equation}
where $\lambdaleft = \lambda + \sum_{\alpha} A(\alpha)\, \alpha$,
$\lambdaright = \lambda + \sum_{\alpha} C(\alpha)\, \alpha$ (both sums
over $\Phi^+$). Moreover, we have equalities
\begin{equation*}
\begin{aligned}
\sum_{\alpha\in \Phi^+} A(\alpha)\, \alpha &=
\sum_{\alpha=\varepsilon_i-\varepsilon_j;\; i<j} A(\alpha)
(\varepsilon_i - \varepsilon_j) \\
&= \ct_L(e_A) - \ct_R(e_A) \\
&= -\ct_L(f_A) + \ct_R(f_A)
\end{aligned}
\end{equation*}
from which it follows that
\begin{equation}
\lambdaleft = \lambda - \ct_R(e_A) + \ct_L(e_A), \quad
\lambdaright = \lambda + \ct_R(f_C) - \ct_L(f_C).
\end{equation}

\begin{lem}
$\ct(e_A 1_\lambda f_C) \preceq \lambda \Leftrightarrow 
\ct_L(1_\lambdaleft e_A f_C) \preceq \lambdaleft  \Leftrightarrow
\ct_R(e_Af_C 1_\lambdaright) \preceq \lambdaright$.
\end{lem}

\begin{proof}
By the definitions we have $\ct(e_A 1_\lambda f_C) = \ct(e_A f_C)$,
with similar equalities for $\ct_R$, $\ct_L$. 
From equation \eqref{content:secondrel} and the above we have the
following equivalences
\begin{align*}
\ct(e_A 1_\lambda f_C) \preceq \lambda
& \Leftrightarrow \ct_R(e_A) + \ct_L(f_C) \preceq \lambda \\
& \Leftrightarrow \ct_L(f_C) \preceq \lambda - \ct_R(e_A) \\
& \Leftrightarrow \ct_L(f_C) + \ct_L(e_A) \preceq \lambda - \ct_R(e_A)
  + \ct_L(e_A) \\
& \Leftrightarrow \ct_L(1_\lambdaleft e_A f_C) \preceq \lambdaleft .
\end{align*}
This establishes the first equivalence of the lemma.
The second equivalence is established by a similar argument.
\end{proof}

By similar reasoning one can obtain similar equivalences
in which the $f_A$ precede the $e_C$. 

It follows from the preceding lemma that the set $\basis_+$ (see
Theorem \ref{thm:idemp:basis}) can be rewritten in either of the forms
\begin{equation}\label{basis:rewrite}
\begin{aligned}
\textstyle\bigcup_\lambdaleft\, 
\{1_\lambdaleft e_A f_C &\mid \ct_L(e_A f_C) \preceq \lambdaleft\}\\
&= \textstyle\bigcup_\lambdaright\, 
\{e_A f_C 1_\lambdaright \mid \ct_R(e_A f_C) \preceq \lambdaright\}
\end{aligned}
\end{equation}
with a similar statement applying to the set $\basis_-$.

The subspace of $U$ spanned by the $H_i$ and $x_\alpha$ is isomorphic
with the Lie algebra $\gl_n$ of $n\times n$ matrices.  The isomorphism
is determined by sending $e_i \to e_{i,i+1}$, $f_i \to
e_{i+1,i}$. Here the notation $e_{ij}$ stands for the matrix with all
entries $0$, except for the $(i,j)$th entry, which is $1$. It follows
from the definition of $x_\alpha$ (see \sect \ref{sec:main}) that the
isomorphism carries $H_i$ to $e_{ii}$ and $x_\alpha$ to $e_{ij}$ when
$\alpha=\varepsilon_i-\varepsilon_j \in \Phi$.  From this one can now
easily verify that (for $\alpha = \varepsilon_i - \varepsilon_j$,
$\beta = \varepsilon_k - \varepsilon_\ell$ for $i\ne j$, $k\ne \ell$)
\begin{equation}\label{gl:relations}
x_\alpha x_\beta - x_\beta x_\alpha = 
\begin{cases} H_\alpha & \text{if $\alpha + \beta = 0$}\\
c_{\alpha,\beta} x_{\alpha+\beta} & \text{if $\alpha+\beta\in \Phi$}\\
0 & \text{otherwise}
\end{cases}
\end{equation}
where $H_\alpha{}=H_i - H_j$ and where for $\alpha+\beta\in \Phi$ we have
$$
c_{\alpha,\beta} = 
\begin{cases}
1 & \text{if $j=k$ and $i\ne \ell$}\\
-1 & \text{if $i=\ell$ and $j\ne k$}.
\end{cases}
$$
The above relations hold in the enveloping algebra $U$, since we use
only the defining relations \eqref{R1}--\eqref{R5} in their
verification. Hence they are valid as well as in the quotient $\B$.

For $r,s \in \N$, $\alpha, \beta \in \Phi$, one deduces the following
commutation relations from the above by induction (compare with
Kostant \cite[Lemma 1]{Kostant}):
\begin{gather} 
\divided{x_\alpha}{r} \divided{x_{-\alpha}}{s} =
\sum_{j=0}^{\min(r,s)} \divided{x_{-\alpha}}{s-j} \binom{H_\alpha-r-s+2j}{j}
\divided{x_\alpha}{r-j}, \tag{\theequation a}\label{commute:xx:0}
\end{gather}
and, for $\alpha + \beta \ne 0$, $\alpha+ \beta \in \Phi$ 
\begin{equation}\label{commute:xx:not0}
\divided{x_\alpha}{r} \divided{x_\beta}{s} =
\divided{x_\beta}{s} \divided{x_\alpha}{r} +
\sum_{j=1}^{\min(r,s)} c_{\alpha,\beta}^j \divided{x_\beta}{s-j}
\divided{x_{\alpha+\beta}}{j} \divided{x_\alpha}{r-j}.\tag{\theequation b}
\end{equation}
For $\alpha + \beta \ne 0$, $\alpha + \beta \notin \Phi$ 
\begin{equation}\label{commute:xx:last}
\divided{x_\alpha}{r} \divided{x_\beta}{s} = 
\divided{x_\beta}{s} \divided{x_\alpha}{r}.\tag{\theequation c}
\end{equation}
From these formulas it follows that we can always interchange the
order of a product of two divided powers of root vectors, up to a
$\Z$-linear combination of terms of strictly lower degree and (right
or left) content.

In the following argument we make use of the fact that
$\dbinom{H_\alpha-t}{s}$ ($s,t \in \Z$, $\alpha\in \Phi$) belongs to
the subring of $U$ generated by the divided powers of all root
vectors, and thus belongs to $\B^0_\Z$.  This follows from \cite[Lemma
1]{Kostant} by an inductive argument.  We will also use the identity
\begin{equation}\label{divided:combine}
\divided{x_\gamma}{a} \divided{x_\gamma}{b} 
= \binom{a+b}{a}\divided{x_\gamma}{a+b} \quad (\gamma\in \Phi). 
\end{equation}
which follows immediately from the definitions.

\begin{prop}\label{prop:Omega:spans}
The sets $\basis_+$ and $\basis_-$ span the algebra $\B_\Z$.
\end{prop}

\begin{proof}  We prove just the claim about $\basis_+$, as the other case is
similar.  We use the second formulation of $\basis_+$ given in
\eqref{basis:rewrite}.  The algebra $\B_\Z$ is spanned by the set of
all products formed from divided powers of root vectors and
idempotents $1_\lambda$.  By Proposition \ref{prop:rootXidemp}, in any
such product, one may always commute the idempotents all the way to
the right.  Hence, $\B_\Z$ is spanned (over $\Z$) by the set of all
monomials of the form
\begin{equation}\label{special}
M = \divided{x_{\gamma_1}}{\psi_1} \cdots
\divided{x_{\gamma_m}}{\psi_m} 1_\mu
\end{equation}
for various $m\in \N$, $\gamma \in \Phi^m$, $\psi \in \N^m$, $\mu\in
\Lambda(n,d)$.  We can assume that $\gamma_i \ne \gamma_{i+1}$ for all
$i\le m-1$, for otherwise the monomial is an integral multiple of a
monomial having that property, by \eqref{divided:combine}.  Call
monomials of such form {\em special}.

Let $\chi=\ct_R(M)$.  We define the {\em deviation} of $M$ by
$\delta{}= \sum_{\chi_i>\mu_i} (\chi_i-\mu_i)$.  Note that if
$$
M^\prime = 
\divided{x_{\gamma^\prime_1}}{\psi^\prime_1} \cdots
\divided{x_{\gamma^\prime_{m'}}}{\psi^\prime_{m'}} 1_\mu
$$
is another special monomial with the same idempotent part $1_\mu$,
then $\ct_R(M) \preceq \ct_R(M') \Rightarrow  \delta(M)\le\delta(M')$
and $\ct_R(M) \preceq \mu \Leftrightarrow \delta(M) = 0$. 

Now we claim that all special monomials $M$ of deviation $0$ lie in
the $\Z$-span of $\basis_+$. We argue by induction on the degree
$r=\sum_i \psi_i$ of $M$. If $r=0$ then $M=1_\mu$, which belongs to
$\basis_+$.  Suppose that $r>0$. If $M\in \basis_+$ then we are
done. Otherwise, we can apply the commutation relations 
\eqref{gl:relations} to reorder the factors in
$M$, obtaining an equality of the form
\begin{align*}
M = \divided{x_{\gamma_1}}{\psi_1} \cdots
\divided{x_{\gamma_m}}{\psi_m} = c\, e_A f_C 
+ \text{lower terms}
\end{align*}
for some $c\in \Z$, $A,C \in \N^{\Phi^+}$.  Here, the lower terms are
integral multiples of Kostant monomials of strictly lower degree and
content.  The constant $c\in \Z$ arises not from the commutation formulas
but when two divided powers are combined via the equality
\eqref{divided:combine}.
Note that, in the lower order terms, whenever a factor of
the form $\dbinom{H_\alpha-t}{s}$ ($s,t\in \Z$) appears, we express
it in terms of a $\Z$-linear combination of $1_\lambda$'s, and then
apply Proposition \ref{prop:rootXidemp} to commute the idempotents as
far to the right as possible. Upon right multiplication of the above
equality by $1_\mu$ we obtain the equality
$$
M 
= c\, e_A f_C 1_\mu  + \text{lower terms}
$$
where the lower terms are integral multiples of terms of the same form
as $M$ (each having a factor $1_\mu$ on the right), but again, of
strictly lower degree and content than that of $M$.  Now, in the above
equality the right content of $e_Af_C1_\mu$ is equal to the right
content of $M$, and thus $e_Af_C1_\mu \in \basis_+$. By induction the
lower terms belong to the $\Z$-span of $\basis_+$. The claim is
proved.

Now we proceed by induction on deviation.  (The above claim forms the
base step in the induction.)  Let $M$ be a special monomial of the
form \eqref{special} of deviation $\delta=\delta(M)>0$. Set
$\chi=\ct_R(M)$.  Since $\chi \not\preceq \mu$, there exists an index
$j$ for which $\chi_j > \mu_j$. Fixing this index $j$, we call
$\beta\in \Phi$ {\em bad} if $\beta=\varepsilon_i - \varepsilon_j$ for
some $i\ne j$ and call $\beta$ {\em good} otherwise.  We extend this
terminology to the factors $\divided{x_\gamma}{\psi_\gamma}$ of $M$.
We define the {\em defect} $D$ of the monomial $M$ by the equality
$D(M) = \sum \psi_\alpha$ where the sum is taken over the set $\{
\divided{x_\alpha}{\psi_\alpha} \}$ of good factors in
$M$ which appear to the right of some bad factor.  Note that $D=0$ if
and only if all the bad factors in $M$ appear as far to the right as
possible.  From Proposition \ref{prop:rootXidemp} it follows
immediately that $M=0$ in $\B$ whenever $D(M)=0$.

Now suppose that $D(M)>0$.  So there exists at least one good factor
$\divided{x_\alpha}{a}$ in $M$ appearing to the right of some bad
factor $\divided{x_\beta}{b}$. We may assume that
$\divided{x_\alpha}{a}$ is the leftmost such good factor in $M$.  It
has one or more bad factors immediately to its left.  We successively
commute $\divided{x_\alpha}{a}$ with each bad factor to its left,
using relations \eqref{gl:relations}.  As before, we express
factors of the form $\dbinom{H_\alpha-t}{s}$ in terms of idempotents
$1_\lambda$, and commute these all the way to the right using
Proposition \ref{prop:rootXidemp}.  By Proposition \ref{prop:Hmult} we
know that such factors act as integral scalars on $1_\mu$. The result
of all this is thus, up to a $\Z$-linear combination $L$ of Kostant
monomials (all involving the same idempotent $1_\mu$ on the right) of
strictly lower right content, an integral multiple of a monomial $M'$
of the same right content as $M$ but of strictly lower defect.  By
induction (on defect) $M'$ lies in the $\Z$-span of $\basis_+$, and,
since the monomials in $L$ have lower deviation than $M$ does, $L$
must also lie within the $\Z$-span of $\basis_+$.  This proves that
$M$ lies in the $\Z$-span of $\basis_+$, and hence that $\basis_+$
spans $\B_\Z$.  The argument for $\basis_-$ is similar (interchange
right and left, $+$ and $-$ in the above argument).
\end{proof}


\section{Proof of Theorems \ref{thm:present:S} 
and \ref{thm:idemp:basis}} \label{sec:pf:thms1&3}

Since $S=S(n,d)$ is the image of the map $U \to \End(V^{\otimes d})$,
it follows immediately from Lemma \ref{lem:minpoly} that the images of
the $H_i$ in the Schur algebra $S$ satisfy the defining relations for
the algebra $\B$, so the surjection $U \to S$ factors through the
algebra $\B$.  In particular, this gives a surjection $\B \to S$.  It
follows that $\dim \B \ge \dim S$.  In order to produce the opposite
inequality, which will complete the proof of Theorem
\ref{thm:present:S}, it is enough to produce a spanning set in $\B$ of
cardinality equal to the dimension of $S$.  We know from Proposition
\ref{prop:Omega:spans} that the sets $\basis_+$, $\basis_-$ span $\B$,
so the proof of Theorem \ref{thm:present:S} is completed by the
following lemma.

\begin{lem} \label{lem:card}
The cardinality of $\basis_+$ and $\basis_-$ is equal to the dimension
of $S=S(n,d)$.
\end{lem}

\begin{proof}
By symmetry, it is enough to prove this for $\basis_+$. It is well
known (see \cite{Green}) that the dimension of $S(n,d)$ is given by
the number of monomials of total degree $d$ in $n^2$ variables.  This
is the same as the number of monomials in $n^2-1$ variables of total
degree not exceeding $d$; in other words, the dimension of $S(n,d)$ is
the same as the cardinality of the set
$$ P{}=
\{ e_A H_B f_C \mid B\in \N^n_1, A,C \in \N^{\Phi^+},
 |A|+|B|+|C|\le d \}.
$$
Thus, to prove the result it suffices to
give a bijective correspondence between $P$ and $\basis_+$. One such
is given by the map
$$
e_A H_B f_C  \to  e_A 1_\lambda f_C 
$$
where $\lambda=(d-|A|-|B|-|C|)\varepsilon_1+B+\ct(e_Af_C)$.  The
inverse map is given by
$$
e_A 1_\lambda f_C  \to  e_A H_B f_C
$$
where $B = \lambda - \ct(e_Af_C) - \lambda_1\varepsilon_1$.  The lemma
is proved.
\end{proof}

It remains to prove Theorem \ref{thm:idemp:basis}.
It follows from the preceding arguments that the quotient map $\B \to
S$ is an isomorphism of algebras, and from Proposition
\ref{prop:Omega:spans} we conclude that $\basis_+$ and $\basis_-$ are bases
for $\B$ (over $\Q$).  Hence these sets are linearly independent over
$\Q$, and thus also over $\Z$.  Thus they are $\Z$-bases for $\B_\Z$.
Carter and Lusztig \cite[Thm.~3.1]{CL} showed that the restriction to
$U_\Z$ of the map $U \to S$ gives a surjection $U_\Z \to
S_\Z{}=S_\Z(n,d)$.  It follows that the restriction map $\B_\Z \to
S_\Z$ is an isomorphism, and that the sets $\basis_+$, $\basis_-$ are
integral bases for the Schur algebra $S_\Z$.  Moreover, the
restriction of the map $U_\Z \to S_\Z$ to $U_\Z(\sl_n)$ is still
surjective, according to \cite[p.~44]{Donkin:SA3}, and thus the image
is generated by all $\divided{e_i}{m}$, $\divided{f_i}{m}$.  This
proves Theorem \ref{thm:idemp:basis}.

\begin{remark}
We conjecture that the set $P$ appearing in the proof of Lemma
\ref{lem:card} is actually another integral basis for $S_\Z(n,d)$. 
More generally, for any fixed $i_0$, ($1\le i_0 \le n$), either of 
the sets
$$
\{ e_A H_B f_C \},  \quad \{ f_A H_B e_C \}
$$
($B\in \N^n_{i_0}, A,C \in \N^{\Phi^+}, |A|+|B|+|C|\le d$)
should be an integral basis of $S_\Z(n,d)$.  This would be a truncated
form of Kostant's well-known basis for $U_\Z$.  The conjecture is
true when $n=2$; see \cite{DG}.
\end{remark}

\newcommand{\XX}{\mathcal{B}}
\section{Proof of Theorem \ref{thm:idemp:present}}
\label{sec:idemp:present}

Let $\XX$ be the $\Q$-algebra given by the generators and
relations of Theorem \ref{thm:idemp:present}.  In $\XX$ we define
elements $H_j$ for $j=1,\dots,n$ by setting $H_j{}=\sum_\lambda
\lambda_j 1_\lambda$, where the sum is carried over all $\lambda\in
\Lambda(n,d)$.

Since the various $1_\lambda$'s commute it follows that the $H_j$'s
must also commute, so relation \eqref{R1} holds for the elements
$H_j$.  Relation \eqref{R2} follows immediately from the defining
relations \eqref{S3} and the definition of the $H_j$.

From the defining relations \eqref{S1} and the definition of
the $H_j$ it follows that
$$
\sum_{j=1}^n H_j = \sum_{j=1}^n \sum_\lambda \lambda_j 1_\lambda
= \sum_\lambda \left(\sum_{j=1}^n \lambda_j\right) 1_\lambda
= \sum_\lambda d\,1_\lambda = d\cdot 1 = d.
$$
This proves that the $H_j$ satisfy relation \eqref{R6}. Moreover,
we also have equalities
$$
1_\lambda \, H_j = H_j \, 1_\lambda = \sum_\mu \mu_j 1_\mu \, 1_\lambda
 = \lambda_j 1_\lambda
$$
for each $\lambda$, $j$, from which we obtain the equalities
\begin{equation}\label{H:choose:b}
\begin{aligned}
H_j(H_j-1)\cdots(H_j-b)
&= \left(\sum_\lambda \lambda_j\, 1_\lambda\right)(H_j-1)\cdots(H_j-b) \\
&= \sum_\lambda \lambda_j(\lambda_j-1)\cdots(\lambda_j-b) \,1_\lambda
\end{aligned}
\end{equation}
for any $b\in \N$.  This is $0$ when $b=d$, since $\lambda$ is a
composition of $d$, and thus each part $\lambda_j$ of $\lambda$ is an
integer in the interval $0,\dots,d$.  This proves that the $H_j$
satisfy relation \eqref{R7}.

We now want to show that the $H_j$ also satisfy relation
\eqref{R3}. For this we will use the defining relations
\eqref{S2}. For convenience, we extend the definition of the symbol
$1_\lambda$ to all $\lambda\in \Z^n$ such that $|\lambda|=\sum
\lambda_i = d$, defining it to have the value $0$ if any part of
$\lambda$ is negative.  With this convention we have
\begin{equation}
H_j e_i = \sum_\lambda \lambda_j\, 1_\lambda e_i  
= \sum_\lambda \lambda_j\, e_i 1_{\lambda-\alpha_i}
\end{equation}
where the sums are taken over all $\lambda\in \Z^n$ such that
$|\lambda|=d$.  Replacing $\lambda-\alpha_i$ by $\mu$ and noting that
$\lambda_j = \mu_j$ if $j \ne i,i+1$, $\lambda_j = \mu_j + 1$ if $j=i$,
and $\lambda_j = \mu_j - 1$ if $j=i+1$ we obtain

\begin{equation}
H_j e_i = \sum_\mu \lambda_j\, e_i 1_\mu 
= e_i H_j + (\delta_{ij}-\delta_{i+1,j})e_i.
\end{equation}
where again the sum is over all $\mu\in \Z^n$ satisfying
$|\mu|=d$.
This proves that the $H_j, e_i$ satisfy relation \eqref{R3}; a
similar argument shows the same for the $H_j, f_i$.

From \eqref{H:choose:b} it follows (upon replacing $b$ by $b-1$ and
dividing by $b!$) that
\begin{equation}
\binom{H_j}{b} = \sum_\lambda \binom{\lambda_j}{b} \, 1_\lambda
\end{equation}
where the sum is over $\Lambda(n,d)$.  It then follows from relations
\eqref{S1} and the above, by a simple calculation, that for any $\mu
\in \Lambda(n,d)$ we have
$$
\prod_{j=1}^n  \binom{H_j}{\mu_j} = 1_\mu
$$
and thus the $H_j$, ($1\le j \le n$) together with the $e_i,f_i$ ($1
\le i \le n-1$) generate the algebra $\XX$.  Since we have proved that
these generators satisfy the defining relations for the algebra $\B$,
it follows that $\XX$ is a homomorphic image of $\B$.

On the other hand, by Proposition \ref{prop:rootXidemp} we know that
the elements $e_i,f_i, 1_\lambda$ of $\B$ satisfy relations
\eqref{S2}. They also satisfy relation \eqref{S1}, clearly, and
relation \eqref{S3}, by Proposition \ref{prop:Hmult}.  Moreover,
$\B$ is generated by the $1_\lambda, e_i, f_i$ since by Proposition
\ref{prop:Hmult}(c) we know that $H_j$ is expressible as a linear
combination of the $1_\lambda$.  Thus $\B$ is a homomorphic image of
$\XX$.  Combining this with the conclusion of the preceding paragraph,
we see that $\XX \cong \B \cong S(n,d)$. Theorem
\ref{thm:idemp:present} is proved.

\newcommand{\BB}{\mathbf{T}}
\section{The algebra $\BB$.} \label{qsec:B}
We turn now to the quantum case.  Fix $n$ and $d$, and set $\S =
\S(n,d)$.  We define an algebra $\BB=\BB(n,d)$ (over $\Q(v)$) by the
generators and relations of Theorem \ref{qthm:present:S}.  Since $\U$
is the algebra on the same generators but subject only to relations
\eqref{Q1}--\eqref{Q5}, we have a surjective quotient map $\U \to
\BB$.  Eventually we shall show that $\BB \simeq \S$, which will prove
Theorem \ref{qthm:present:S}.

The $q$-analogue of Lemma \ref{lem:minpoly} is the following. The
proof is similar to the proof in the classical case, except that 
it is multiplicative where the classical argument is additive. 

\begin{lem}\label{qlem:minpoly}
Under the representation $\U \to \End(\V^{\otimes d})$ the images of
the $K_i$ satisfy the relations \eqref{Q6} and \eqref{Q7}.  Moreover,
the relation \eqref{Q7} is the minimal polynomial of $K_i$ in
$\End(\V^{\otimes d})$.
\end{lem}

Since $\S = \S(n,d)$ is the image of the representation $\U \to
\End(\V^{\otimes d})$, we have a surjection $\U \to \S$.  From the
lemma it follows that the surjection $\U \to \S$ factors through
$\BB$. Because $\BB$, $\S$ are homomorphic images of $\U$, any relations
between generators in $\U$ will carry over to $\BB$, $\S$.  We do not
distinguish notationally between the generators or root vectors for
$\U$, $\BB$, or $\S$.

Rosso \cite{Rosso} has shown that multiplication defines a
$\Q(v)$-linear isomorphism 
$\U^- \otimes \U^0 \otimes \U^+ \xrightarrow{\approx} \U$,
where $\U^+$ (resp., $\U^-$) is the subalgebra of $\U$ generated by
the $E_i$ (resp., $F_i$), and $\U^0$ is the subalgebra of $\U$
generated by all $K_i, K_i^{-1}$. 
It follows that $\U = \U^-\U^0\U^+$.
From this we obtain a similar triangular decomposition of $\BB$:
\begin{equation}\label{BB:TD}
\BB = \BB^- \BB^0 \BB^+
\end{equation}
where $\BB^+, \BB^-, \BB^0$ are defined to be the images of $\U^+,
\U^-, \U^0$ under the quotient mapping $\U \to \BB$.

There are similar factorizations over $\A$. Setting
$\U_\A^+$, $\U_\A^-$, $\U_\A^0$ to be, respectively, the intersection
of $\U^+$, $\U^-$, $\U^0$ with Lusztig's $\A$-form $\U_\A$ (the
$\A$-subalgebra of $\U$ generated by the $\divided{E_i}{m}$,
$\divided{F_i}{m}$, $K_i^{\pm1}, \sqbinom{K_i}{m}$), Du
\cite[\sect2]{Du} has shown (using results of Lusztig) that $\U_\A =
\U_\A^- \U_\A^0 \U_\A^+$; thus
\begin{equation}\label{BBA:TD}
\BB_\A = \BB_\A^- \BB_\A^0 \BB_\A^+
\end{equation}
where the various subalgebras are defined in the obvious way as
appropriate homomorphic images of $\U_\A^+$, $\U_\A^-$, $\U_\A^0$.

Since $\U_\A^+$ (resp., $\U_\A^-$) is the $\A$-subalgebra of $\U$
generated by the $\divided{E_\alpha}{m}$ (resp.,
$\divided{F_\alpha}{m}$) for $\alpha\in \Phi^+$ and $m\in \N$, the
same statement applies to $\BB_\A^+$ (resp., $\BB_\A^-$) in relation
to $\BB$. Moreover, $\U_\A^0$ is the $\A$-subalgebra of $\U$ generated
by the $K_i^{\pm1}$, $\sqbinom{K_i}{m}$ for $1\le i\le n$, $m\in
\N$, so $\BB_\A^0$ is the $\A$-subalgebra of $\BB$ generated by the
same elements.

Now we determine the structure of the algebra $\BB^0$.  As we shall
see, the structure turns out to be essentially the same as that in the
classical case.  Consider first the algebra $\U^0$, which may be
identified with the commutative polynomial algebra $\Q(v)[K_1^{\pm 1},
\ldots , K_n^{\pm 1}]$. We define $\S^0$ to be the image of $\U^0$
under the quotient map $\U \to \S$. As in the classical case, this map
factors through the algebra $\BB^0$.

\begin{prop} \label{prop:qidempotent}
Define an algebra $\BB^\prime = \U/I^0$ where $I^0$ is the ideal in
$\U^0$ generated by the elements $(K_i-1)(K_i-v)\cdots (K_i-v^d)$
($1\le i \le n$) and $K_1K_2\cdots K_n -v^d$.

\par\noindent(a)  We have an algebra isomorphism $\BB^\prime \cong \BB^0$. 

\par\noindent(b) The set $\{ 1_\lambda \mid \lambda \in \Lambda(n,d)\}$
is a $\Q(v)$-basis for $\BB^0$ and an $\A$-basis for $\BB_\A^0$;
moreover, this set is a set of pairwise orthogonal idempotents which
add up to $1$.

\par\noindent(c) $K_B = 0$ for any $B\in \N^n$ such that $|B|>d$.
\end{prop}

\begin{proof} The argument is similar to the proof of Proposition
\ref{prop:idempotent}.  Consider first the algebra
$\widetilde{\BB}^\prime$ defined to be the quotient of $\U^0$ by the ideal
generated only by the $p(K_i){}=(K_i-1)(K_i-v)\cdots (K_i-v^d)$
for $i=1,\ldots,n$. Since $p(K_i)$ is a non-constant polynomial
with non-zero constant term, each $K_i^{-1}$ already lies in the
ring $\widetilde{\BB}^\prime$.  Thus it follows that
$$
\widetilde{\BB}^\prime \cong \Q(v)[K_1]/(p(K_1))
\otimes \cdots \otimes \Q(v)[K_n]/(p(K_n)).
$$
From the Chinese Remainder Theorem we obtain, as in the proof of
\ref{prop:idempotent}, an isomorphism (products denote direct
products)
$$
\widetilde{\BB}^\prime  \cong \prod_{0\leq \mu_1,\ldots \mu_n \leq d}
\Q(v)[K_1,\ldots , K_n]/(K_1-v^{\mu_1},\ldots , K_n-v^{\mu_n}).
$$
Since $\BB^\prime =\widetilde{\BB}^\prime/(K_1\cdots K_n-v^d)$, we
deduce that
$$
\BB^\prime \cong  \prod_{\mu \in \Lambda(n,d)}
\Q(v)[K_1\ldots ,K_n]/(K_1-v^{\mu_1},\ldots , K_n-v^{\mu_n}).
$$
The preceding isomorphism, which we denote by $\phi$, is given by the
map sending $f(K_1,\ldots K_n)$ to the vector $(f(v^{\mu_1},\ldots
,v^{\mu_n}))_{\mu \in \Lambda(n,d)}$.  In particular, if $\lambda \in
\Lambda(n,d)$, then
$$
\phi (1_\lambda)= \left(\sqbinom{\mu_1}{\lambda_1} \cdots \sqbinom
{\mu_n}{\lambda _n}\right)_{\mu \in \Lambda(n,d)} = (\delta_{\lambda
\mu})_{\mu\in \Lambda(n,d)}.
$$
Thus the various $1_{\lambda}$ are orthogonal idempotents whose sum is
the identity.

Now let $I$ be the ideal in $\U$ generated by elements
$(K_i-1)(K_i-v)\cdots(K_i-v^d)\ \ (1\le i \le n)$ and $K_1 \cdots
K_n-v^d$.  Then by definition $\BB \cong \U/I$.  The canonical quotient
map $\U \to \U/I$ induces, upon restriction to $\U^0$, a map $\U^0 \to
\U/I$.  The image of this map is $\BB^0$ and its kernel is $\U^0 \cap I$,
so $\BB^0 \cong \U^0/(\U^0 \cap I)$.  Clearly $I^0 \subset \U^0 \cap I$.
Thus we obtain the following sequence of algebra surjections
\begin{equation}
\begin{CD}
\BB^\prime = \U^0/I^0 @>>> \U^0/(\U^0 \cap I) @>\approx>> \BB^0 @>>> \S^0
\end{CD}
\end{equation}
where the last map is obtained by Lemma \ref{qlem:minpoly}, and the
middle map is actually an isomorphism. By the above we see that the
dimension of $\BB^\prime$ is the cardinality of the set
$\Lambda(n,d)$. This is the same as $\dim \S^0$, so all the
surjections above are algebra isomorphisms.  This proves parts (a) and
(b).

Part (c) is proved in exactly the same way as part (c) of Proposition
\ref{prop:idempotent}.
\end{proof}

\begin{prop} \label{qprop:Kmult}
Suppose $1 \le i \le n$, $c\in \Z$, $t\in \N$, $\lambda \in
\Lambda(n,d)$, and $B \in \N^n$.  Then we have the following
identities in the algebra $\BB^0$:
\begin{gather}
K_i^{\pm 1} 1_\lambda = v^{\pm \lambda_i} 1_\lambda; \quad 
\sqbinom{K_i;c}{t} 1_\lambda = \sqbinom{\lambda_i+c}{t} 1_\lambda\tag{a}\\
K_B \, 1_\lambda = \lambda_B\,1_\lambda, \quad
\text{where $\lambda_B {}= \prod_i \sqbinom{\lambda_i}{B_i}$}\tag{b}\\
K_B = \sum_{\lambda} \lambda_B \, 1_\lambda,\tag{c}
\end{gather}
where the sum in part (c) is carried out over all $\lambda \in
\Lambda(n,d)$.
\end{prop}

\begin{proof}
We prove part (a). Apply the isomorphism $\phi$ from the proof of
Proposition \ref{prop:qidempotent}. We have $\phi(K_i^{\pm 1}
1_\lambda) = (v^{\pm \mu_i} \delta_{\lambda\mu})_\mu = v^{\pm
\lambda_i} \phi(1_\lambda)$.  It follows that $K_i^{\pm 1} 1_\lambda =
v^{\pm \lambda_i} 1_\lambda$.  Similarly, we have 
\begin{align*}
\phi(\sqbinom{K_i;c}{t} 1_\lambda) &= \phi\left( \left(\prod_{s=1}^t
\frac{K_iv^{c-s+1} - K_i^{-1}v^{-c+s-1}}{v^s-v^{-s}} \right) 1_\lambda
\right) \\ 
&= \left( \prod_{s=1}^t \frac{v^{\mu_i}v^{c-s+1} -
v^{-\mu_i}v^{-c+s-1}}{v^s-v^{-s}} \delta_{\lambda\mu} \right)_\mu \\
&= \left( \sqbinom{\mu_i+c}{t} \delta_{\lambda\mu} \right)_\mu \\
&= \sqbinom{\lambda_i+c}{t} \phi(1_\lambda)
\end{align*}
which proves the second equality in (a). Note that the equality
\begin{equation}
\sqbinom{K_i}{t} 1_\lambda = \sqbinom{\lambda_i}{t} 1_\lambda
\end{equation}
is the case $c=0$ of the above.

The rest of the proof is similar to the proof of Proposition
\ref{prop:Hmult}.
\end{proof}

\noindent
By essentially the same argument as in the classical case we
obtain

\begin{cor}
For any fixed choice of $i_0$ ($1\le i_0 \le n$) the set
$\{ K_B \mid B \in \N^n_{i_0}, |B| \le d \}$ is a $\Q(v)$-basis for
$\BB^0$ and an $\A$-basis for $\BB_\A^0$.
\end{cor}

We also have the following exact analogue of Proposition
\ref{prop:rootXidemp}. 

\begin{prop} \label{qprop:rootXidemp}
For any $\alpha \in \Phi$ and any $\lambda \in \Lambda(n,d)$
we have the commutation formulas
$$
X_\alpha 1_\lambda =
\begin{cases}
1_{\lambda+\alpha} X_\alpha & \text{if $\lambda+\alpha \in \Lambda(n,d)$}\\
0 & \text{otherwise}
\end{cases}
$$
and similarly
$$
1_\lambda X_\alpha =
\begin{cases}
X_\alpha 1_{\lambda-\alpha} & \text{if $\lambda-\alpha \in \Lambda(n,d)$}\\
0 & \text{otherwise}.
\end{cases}
$$
\end{prop}

\begin{proof}
We will need the following identities (see \cite[\sect2.3 (g3),
(g4)]{Lusztig}):
\begin{gather}
\sqbinom{K_i;0}{t} \sqbinom{K_i;-t}{t'} 
= \sqbinom{t+t'}{t} \sqbinom{K_i;0}{t+t'} \quad(t,t' \in \N)\label{ident:a}\\
\sqbinom{K_i;c+1}{t} = v^t \sqbinom{K_i;c}{t} 
+ v^{t-c-1}K_i^{-1} \sqbinom{K_i;c}{t-1} \quad(t \ge 1). \label{ident:b}
\end{gather}
From the defining relation \eqref{Q3} and the definition (see
\sect\ref{qsec:main}) of the root vector $X_\alpha$ one proves that
\begin{equation}\label{qcommute:KX}
K_l X_\alpha = v^{(\varepsilon_l,\alpha)} X_\alpha K_l
\qquad(\alpha\in \Phi).
\end{equation}
From this it follows by a simple calculation that
\begin{equation}\label{qcommute:XKsemi}
\sqbinom{K_l;c}{t} X_\alpha 
= X_\alpha \sqbinom{K_l;c+(\varepsilon_l,\alpha)}{t}
\qquad(\alpha\in \Phi, c\in \Z, t\in \N).
\end{equation}
From this last equality it follows that for $\lambda \in
\Lambda(n,d)$, $\alpha = \varepsilon_i - \varepsilon_j \in \Phi$
($i\ne j$) we have
$$
X_\alpha 1_\lambda = \sqbinom{K_i;-1}{\lambda_i}\sqbinom{K_j;1}{\lambda_j}
\left(\prod_{l\ne i,j} \sqbinom{K_l}{\lambda_l} \right) X_\alpha.
$$
Multiply both sides of the preceding equality by $\sqbinom{K_i}{1} =
\sqbinom{K_i;0}{1}$ and use \eqref{ident:a} to simplify the right-hand-side
and use \eqref{qcommute:XKsemi} to simplify the left-hand-side. The result is
the equality
$$
X_\alpha \sqbinom{K_i;1}{1} 1_\lambda 
= \sqbinom{\lambda_i+1}{1} \sqbinom{K_i}{\lambda_i+1} \sqbinom{K_j;1}{\lambda_j}
\left(\prod_{l\ne i,j} \sqbinom{K_l}{\lambda_l} \right) X_\alpha.
$$
By Proposition \ref{qprop:Kmult}(a) the left-hand-side of this equality
simplifies to give the equality
$$
\sqbinom{\lambda_i+1}{1} X_\alpha 1_\lambda 
= \sqbinom{\lambda_i+1}{1} \sqbinom{K_i}{\lambda_i+1} \sqbinom{K_j;1}{\lambda_j}
\left(\prod_{l\ne i,j} \sqbinom{K_l}{\lambda_l} \right) X_\alpha
$$
and after cancelling the common scalar factor and expanding by means of
\eqref{ident:b} (assuming $\lambda_j > 0$) this becomes
$$
X_\alpha 1_\lambda 
= \sqbinom{K_i}{\lambda_i+1} \left(v^{\lambda_j} \sqbinom{K_j}{\lambda_j} 
+ v^{\lambda_j-1} K_j^{-1}\sqbinom{K_j}{\lambda_j-1} \right) 
\left(\prod_{l\ne i,j} \sqbinom{K_l}{\lambda_l} \right) X_\alpha.
$$
In case $\lambda_j > 0$ after multiplying through in the above
expression and applying Proposition \ref{prop:qidempotent}(c) we see
that the first summand must be zero, so the above equality in that
case simplifies to $X_\alpha 1_\lambda = v^{\lambda_j-1}K_j^{-1}
1_{\lambda+\alpha} X_\alpha$. Now by Proposition \ref{qprop:Kmult}(a)
$K_j^{-1}$ acts on $1_{\lambda+\alpha}$ as $v^{-(\lambda_j-1)}$. Thus
we obtain the equality in the first part of the proposition in the
case $\lambda_j > 0$.

If $\lambda_j = 0$ then one sees easily by Proposition
\ref{prop:qidempotent}(c) that the right-hand-side vanishes. The first
part of the proposition is proved. The proof of the other part is
similar.
\end{proof}

\begin{cor}\label{qcor:root:nil}
We have in $\BB$ the equalities $X_\alpha^{d+1} = 0$ for all $\alpha\in
\Phi$.
\end{cor}

Again, the proof is the same as in the classical case.

\section{Quantum straightening}

Similar to the classical case, we define left and right content
$\ct_L$, $\ct_R$ on Kostant monomials by
\begin{equation}\label{qcontent:firstrel}
\begin{gathered}
\ct_R(\divided{X_\alpha}{m}){}= m\,\varepsilon_j, \quad
\ct_L(\divided{X_\alpha}{m}){}= m\,\varepsilon_i, \\
\ct_R(\sqbinom{K_i}{m} ) = \ct_L(\sqbinom{K_i}{m}) 
= \ct_R(K_i)=\ct_L(K_i){}= 0
\end{gathered}
\end{equation}
where $\alpha = \varepsilon_i - \varepsilon_j$, and again using the
rule $\ct_L(XY)=\ct_L(X)+\ct_L(Y)$ (similarly for $\ct_R$) whenever
$X$ and $Y$ are Kostant monomials.

The $q$-analogues of \eqref{content:firstrel}--\eqref{basis:rewrite}
hold, with the same argument as in the classical case. In particular,
we can work with the description of $\qbasis_+$ in which
the idempotent appears on the right.

Now we write $X_{ij}$ for the root vector $X_\alpha$ when
$\alpha=\varepsilon_i - \varepsilon_j$ and write $K_\alpha = K_{ij} =
K_i K_j^{-1}$ ($1\le i\ne j \le n$). Note that $E_i = X_{i,i+1}$, $F_i
= X_{i+1,i}$.  We assume that $i<j$ and $k<l$. By Xi
\cite[\sect5.6]{Xi} we have the following root vector commutation
formulas, listed below in formulas \eqref{Uq:rel++} and
\eqref{Uq:rel+-} (some of which appeared already in \cite{Lusztig} and
\cite{Rosso}). Note that our notation differs from Xi's: he writes
$E_{i,j-1}$ (resp., $F_{i,j-1}$, $K_{i,j-1}$) where we write $X_{ij}$
(resp., $X_{ji}$, $K_{ij}$).
\begin{equation}\label{Uq:rel++}
X_{ij} X_{kl} =
\begin{cases}
X_{kl} X_{ij} &\hspace{-0.8in} (j < k \text{ or } k < i < j < l) \\
v^{-1} X_{kl} X_{ij} & \hspace{-0.8in}(i = k < j < l 
\text{ or } i < k < j = l) \\
v X_{il} + v X_{kl} X_{ij} &\hspace{-0.4in} (j = k) \\
X_{kl} X_{ij} + (v^{-1}-v) X_{il} X_{kj} & \quad(i < k < j < l) 
\end{cases}
\end{equation}
which, as Xi ({\em loc.\ cit.}) proves, lead to the following
\begin{equation}
\divided{X_{ij}}{M} \divided{X_{kl}}{N} 
= \divided{X_{kl}}{N} \divided{X_{ij}}{M} 
\qquad (j < k \text{ or } k < i < j < l) \tag{\theequation a}
\end{equation}
\begin{equation}
\begin{split}
\divided{X_{ij}}{M} \divided{X_{kl}}{N} 
= v^{-MN} \divided{X_{kl}}{N}\divided{X_{ij}}{M}& \\
(i &= k < j < l \text{ or } i < k < j = l) 
\end{split}\tag{\theequation b}
\end{equation}
\begin{equation}
\begin{split}
\divided{X_{ij}}{M} \divided{X_{kl}}{N} 
= \sum_{t=0}^{\min(M,N)} v^{(M-t)(N-t)+t} \divided{X_{kl}}{N-t} 
\divided{X_{il}}{t} \divided{X_{ij}}{M-t}& \\
(j &= k) 
\end{split}\tag{\theequation c}
\end{equation}
\begin{equation}\label{badcase:one}
\begin{split}
\divided{X_{ij}}{M} \divided{X_{kl}}{N} 
= \sum_{t=0}^{\min(M,N)} v^{-t(t-1)/2} (v^{-1}-v)^t [t]! 
\divided{X_{kj}}{t} \divided{X_{kl}}{N-t} &\divided{X_{ij}}{M-t} 
\divided{X_{il}}{t} \\
&(i < k < j < l). 
\end{split}\tag{\theequation d} 
\end{equation}
We note that these formulas lead to others, of the same form as those
given above, except for signs and (in some cases) a scalar factor of some
integral power of $v$. The new formulas are obtained from the ones
listed above by solving for the term $\divided{X_{kl}}{N}
\divided{X_{ij}}{M}$ and then interchanging $(i,j)$ and $(k,l)$.

We also note that the formulas one obtains in this way, together with
the ones already listed, exhaust the possibilities for a commutation
of two root vectors labelled by positive roots.  Indeed, given two
finite intervals $[i,j]$, $[k,l]$ one of the following mutually
exclusive possibilities must apply: 1) the intervals are disjoint; 2)
one interval is properly included in the other without shared
endpoints; 3) one interval is properly included in the other with the
two intervals sharing a common endpoint; 4) the intervals coincide; 5)
the intervals meet at a single point; 6) the intervals properly
overlap.  Looking at the formulas listed above, we see that case (a)
covers possibilities 1) and 2), case (b) covers possibility 3), case
(c) is possibility 5), and case (d) is possibility 6).  We do not need
a formula for possibility 4).

There are five cases listed in \cite[\sect5.6 (d0)]{Xi} and only four
cases listed above, but the second of the five cases in ({\em loc.\
cit.}) is superfluous, as we have just proved.

There are similar commutation formulas for products of negative root
vectors. One way to get them is by applying the isomorphism $\Omega:
\U \to \U^{\text{opp}}$ of \cite[1.2]{Lusztig}, defined by $\Omega E_i
= F_i$, $\Omega F_i = E_i$, $\Omega K_i = K_i^{-1}$, $\Omega v =
v^{-1}$, to the positive root vector commutation formulas discussed
above.

Xi ({\em loc.\ cit.}) also gives the following, which express the
commutation between positive and negative root vectors.  Still under
the assumption $i<j$ and $k<l$ we have
\begin{equation}\label{Uq:rel+-}
X_{ij} X_{lk} = 
\begin{cases}
X_{lk} X_{ij}  &\hspace{-0.9in} (j \le k \text{ or } k < i < j < l) \\
X_{lk} X_{ij} + v^{-1} K_{kj}^{-1} X_{ik} & (i < k < j = l) \\
X_{lk} X_{ij} - X_{lj} K_{ij}^{-1} & (i = k < j < l) \\
X_{lk} X_{ij} + \sqbinom{K_{ij};0}{1} & (i = k, j = l) \\
X_{lk} X_{ij} + v^{-1}(v-v^{-1}) X_{lj} K_{kj}^{-1} X_{ik} & (i < k < j < l) 
\end{cases}
\end{equation}
which, as Xi proves, lead to the following
\begin{equation}\label{Uq:rel+-a}
\divided{X_{ij}}{M} \divided{X_{lk}}{N} 
= \divided{X_{lk}}{N} \divided{X_{ij}}{M} 
\qquad (j \le k \text{ or } k < i < j < l) \tag{\theequation a}
\end{equation}
\begin{equation}
\begin{split}
\divided{X_{ij}}{M} \divided{X_{lk}}{N} 
= \sum_{t=0}^{\min(M,N)} v^{t(N-t-1)} \divided{X_{jk}}{N-t} K_{kj}^{-t} 
\divided{X_{ij}}{M-t} \divided{X_{ik}}{t}& \\ 
(i < k &< j = l) 
\end{split}\tag{\theequation b}
\end{equation}
\begin{equation}
\begin{split}
\divided{X_{ij}}{M} \divided{X_{lk}}{N} 
= \sum_{t=0}^{\min(M,N)} (-1)^t v^{t(M-t)} \divided{X_{lj}}{t} 
\divided{X_{lk}}{N-t} K_{ij}^{-t} \divided{X_{ij}}{M-t} \\
\qquad (i = k < j < l) 
\end{split}\tag{\theequation c}
\end{equation}
\begin{equation}
\divided{X_{ij}}{M} \divided{X_{ji}}{N} 
= \sum_{t=0}^{\min(M,N)} \divided{X_{ji}}{N-t}
\sqbinom{K_{ij}; 2t-M-N}{t} \divided{X_{ij}}{M-t} \tag{\theequation d}
\end{equation}
\begin{equation}\label{badcase:two}
\begin{split}
\divided{X_{ij}}{M} \divided{X_{lk}}{N} 
&= \sum_{t=0}^{\min(M,N)} \xi(v,t)  
\divided{X_{lk}}{N-t} \divided{X_{lj}}{t} K_{kj}^{-t} 
\divided{X_{ij}}{M-t} \divided{X_{ik}}{t}\\ 
\text{where } \xi(v,t) &= v^{-t(2N+t-1)/2} (v-v^{-1})^t [t]!\qquad  
(i < k < j < l) .
\end{split}\tag{\theequation e}
\end{equation}
As before, these formulas lead to others, of a similar form as those
given above. The new formulas in this case are obtained by first
applying $\Omega$ to get a formula in $\U^{\text{opp}}$, then
switching order of factors to get a formula in $\U$, and finally
interchanging $(i,j)$ with $(k,l)$.

The formulas one obtains in this way, together with
the ones already listed, exhaust the possibilities for a commutation
involving a positive root vector followed by negative root vector.  
For the possibilities 1) -- 6) listed earlier for two finite intervals
$[i,j]$, $[k,l]$, we see that case (a) covers possibilities 1), 5),
and 2), cases (b), (c) cover possibility 3), case (d) is possibility
4), and case (e) is possibility 6). 

There are similar commutation formulas for the case of a negative root
vector followed by a positive root vector. They are obtained from the
ones already described, by simply solving for the term
$\divided{X_{lk}}{N} \divided{X_{ij}}{M}$. The new formulas will be of
the same form, except for signs.

Thus we see that the formulas listed in \eqref{Uq:rel++} and
\eqref{Uq:rel+-}, together with formulas easily derivable from them,
give the commutation in $\U$ between $q$-divided powers of any two
root vectors. 

\medskip

The most important feature of these commutation formulas, for our
purposes, is that every product of the form $\divided{X_\alpha}{r}
\divided{X_\beta}{s}$ may be expressed as a scalar multiple (by an
element of $\A$) of the product in the opposite order, modulo an
$\A$-linear combination of monomials of degree and (right or left)
content which is no greater that that of the original product.  This
differs from the classical case, where the extra terms in any such
commutation always have strictly lower degree and content.

A central result of \cite{Lusztig} is a $q$-analogue of Kostant's
basis for $U_\Z$.  The description of this basis uses the ``box
ordering'' on $\Phi^+$ defined as follows: if $\alpha =\ep_i -\ep_j$
and $\be = \ep_r - \ep_s$, then $\al \succ \be$ if either $s>j$ or
$s=j$ and $r>i$. For $B\in \N^{n-1}$ and $\delta \in \{0,1\}^{n-1}$
set
$$
\KK _{\delta,B} = K_{\alpha_1}^{\delta_1}\cdots
K_{\alpha_{n-1}}^{\delta_{n-1}}\sqbinom{K_{\alpha_1};0}{B_1}\cdots
\sqbinom{K_{\alpha_{n-1}};0}{B_{n-1}}.
$$
Then Lusztig \cite[Thm.~4.5]{Lusztig} proves that the set of all
elements of the form
\begin{equation}\label{Lusztig:basis}
 F_A  \KK _{\delta,B} E_C  \quad (A,C \in \N^{\Phi^+}, B\in \N^{n-1}, 
 \delta \in \{0,1\}^{n-1} ) 
\end{equation}
is an $\A$-basis of $\U_\A(\sl_n)$, provided the products in $E_C$ are
taken in the box order on $\Phi^+$ and the products in $F_A$ are taken
in the reverse box order.  In \cite[Theorem 2.4]{Xi2} it is shown that
the box order is not necessary as one may take any ordering on
$\Phi^+$ when forming $F_A$ and $E_C$.  By applying the involution
$\omega$ (see \cite[3.1.3]{Lusztig:book}) which interchanges the $E_i$
and $F_i$ we obtain another such basis, consisting of all elements of
the form
\begin{equation}\label{Lusztig:altbasis}
 E_A  \KK _{\delta,B} F_C  \quad (A,C \in \N^{\Phi^+}, B\in \N^{n-1}, 
 \delta \in \{0,1\}^{n-1} ) 
\end{equation}

We note that elements $\sqbinom{K_\alpha;c}{t}$ appear in some of
Xi's commutation formulas. 
In the following argument we use the fact that any $\sqbinom{K_\alpha;
c}{t}$ ($\alpha\in \Phi$, $c\in \Z$, $t\in \N$) belongs to $\U^0_\A$
and thus belongs to $\BB^0_\A$.  This was proved by Lusztig 
\cite[4.5]{Lusztig:88} in case $\alpha$ is simple. Lusztig's argument
extends immediately to general $\alpha\in \Phi$ since one has a
version of \cite[(4.3.1)]{Lusztig:88} for any $\alpha$. 

We shall also need the relation
\begin{equation}\label{qdivided:combine}
\divided{X_\gamma}{a} \divided{X_\gamma}{b} 
= \sqbinom{a+b}{a}\divided{X_\gamma}{a+b} \quad (\gamma\in \Phi). 
\end{equation}
which follows immediately from the definitions.

\begin{prop} \label{prop:Omega:qspans}
The sets $\qbasis_+$ and $\qbasis_-$ span the algebra $\BB_\A$.
\end{prop}
\begin{proof} We just prove the statement for $\qbasis_+$.
As in the classical case, the algebra $\BB_\A$ is
spanned by the set of all ``special'' monomials of the form
\begin{equation}\label{qspecial}
M = X_{\gamma_1}^{(\psi_1)} \cdots X_{\gamma_m}^{(\psi_m)} \, 1_{\mu}
\end{equation}
for $m\in \N$, $\gamma\in \Phi^m$, $\psi\in\N^m$.  By
\eqref{qdivided:combine} we may assume that $\gamma_i \ne
\gamma_{i+1}$ for all $i \le m-1$.  Set $\chi =\chi_R(M)$ and define
the deviation of $M$ by $\delta{}=\sum _{\chi_i > \mu_i}(\chi_i -
\mu_i)$.  We claim that the set of all
special monomials of deviation $0$ lie within the $\A$-span of
$\qbasis_+$.  However, the proof of this claim is different in the
present case.  By \eqref{Lusztig:altbasis} we can express the product
$X_{\gamma_1}^{(\psi_1)} \cdots X_{\gamma_m}^{(\psi_m)}$ (which lies
in the subalgebra $\U(\sl_n)$) in terms of an $\A$-linear combination
of terms of the form $E_A \KK _{\delta,B} F_C$.  In this expression,
each of the terms has right content not exceeding that of $M$, because
commutation does not increase degree or right content.  We express
each $\KK _{\delta,B}$ in terms of an $\A$-linear combination of
idempotents $1_\lambda$, and use Proposition \ref{qprop:rootXidemp} to
commute the idempotents all the way to the right.  Upon right
multiplication of the resulting expression by $1_\mu$ we obtain an
equality of the form
$$
M = \sum_{A,C} a_{A,C}\,E_A F_C 1_\mu  \quad (a_{A,C}\in \A)
$$
in which the right content of each term on the
right-hand-side is no greater than that of $M$. It follows that each
term in the right-hand-side of the above equality has deviation $0$,
and thus lies in $\qbasis_+$.  This proves the claim.

Now we proceed by induction on deviation, with the above claim forming
the base step in the induction.  Let $M$ be a special monomial of the
form \eqref{qspecial} of deviation $\delta=\delta(M)>0$. Set
$\chi=\ct_R(M)$.  Since $\chi \not\preceq \mu$, there is an index $j$
for which $\chi_j > \mu_j$. Fixing this index $j$, we call $\beta\in
\Phi$ {\em bad} if $\beta=\varepsilon_i - \varepsilon_j$ for some
$i\ne j$ and call $\beta$ {\em good} otherwise.  We extend this
terminology to the factors $\divided{X_\gamma}{\psi_\gamma}$ of $M$.
We define the {\em defect} $D$ of the monomial $M$ by the equality
$D(M) = \sum \psi_\alpha$ where the sum is over the set $\{
\divided{X_\alpha}{\psi_\alpha} \}$ of good factors in $M$ appearing
to the right of some bad factor. The defect of $M$ is $0$ if and only
if all the bad factors in $M$ appear as far to the right as possible.
As in the classical case we have by Proposition \ref{qprop:rootXidemp}
that $M=0$ whenever $D(M)=0$.

Now suppose that $D(M)>0$.  Then there exists at least one good factor
$\divided{X_\alpha}{a}$ appearing to the right of some bad factor.  We
may assume that $\divided{X_\alpha}{a}$ is the leftmost such good
factor in $M$.  It has one or more bad factors immediately to its
left.  If there is just one bad factor $\divided{X_\beta}{b}$ to the
left of $\divided{X_\alpha}{a}$ we apply relations \eqref{Uq:rel++},
\eqref{Uq:rel+-} to the product
$\divided{X_\beta}{b}\divided{X_\alpha}{a}$. In all cases except
\eqref{badcase:one}, \eqref{badcase:two} the argument is similar to
the classical case, and we obtain a multiple (by an element of $\A$)
of a Kostant monomial of the same deviation as $M$ but of strictly
lower defect, modulo an $\A$-linear combination of monomials of
strictly lower deviation. In cases \eqref{badcase:one},
\eqref{badcase:two} (as written) we obtain an $\A$-linear combination
of monomials of the same defect and deviation as $M$. But in relation
\eqref{badcase:one} the factors $\divided{X_{ij}}{M-t}$,
$\divided{X_{il}}{t}$ commute (up to a power of $v$), and so do the
factors $\divided{X_{kj}}{t}$, $\divided{X_{kl}}{N-t}$.  After
interchanging those pairs of factors, the factors
$\divided{X_{kj}}{t}$, $\divided{X_{il}}{t}$ will be adjacent, and
they commute.  Thus we see that the right-hand-side of formula
\eqref{badcase:one} can be rewritten in such a way that the bad
factors $\divided{X_{kj}}{t}$, $\divided{X_{ij}}{M-t}$ appear on the
right, and thus all the monomials we obtain after applying
\eqref{badcase:one} to $M$ are of strictly lower defect than that of
$M$.  In case \eqref{badcase:two} similar commutation applies to
obtain the same result: here $\divided{X_{ij}}{M-t}$,
$\divided{X_{ik}}{t}$ commute up to a power of $v$, then
$\divided{X_{ik}}{t}$, $\divided{X_{lj}}{t}$ commute by case
\eqref{Uq:rel+-a} since $k-1 < j$.  (One needs also in this case to
express the factor $K_{kj}^{-t}$ as an $\A$-linear combination of
idempotents, and then commute those all the way to the right.)

Now suppose there is more than one bad factor to the left of
$\divided{X_\alpha}{a}$.  Let $\divided{X_\beta}{b}$ be the rightmost
such bad factor. We again apply relations \eqref{Uq:rel++},
\eqref{Uq:rel+-} to the product
$\divided{X_\beta}{b}\divided{X_\alpha}{a}$ and repeat the argument
given above. The result in the cases \eqref{badcase:one},
\eqref{badcase:two} is an $\A$-linear combination of Kostant monomials
of the same defect and deviation as that of $M$, but with fewer bad
factors to the left of the good factor $\divided{X_\alpha}{a}$. By
induction on the number of such factors, we are done.

This proves that $M$ lies in the
$\A$-span of $\qbasis_+$, and hence that $\qbasis_+$ spans $\BB_\A$.
The argument for $\qbasis_-$ is similar.
\end{proof}


\section{Proof of Theorems \ref{qthm:present:S} and
\ref{qthm:idemp:basis}.}

Since $\S$ is the image of the representation $\U \to \End(\V^{\otimes
d})$, it follows from Lemma \ref{qlem:minpoly} that the images of the
$K_i$ in the quantum Schur algebra $\S$ satisfy the defining relations
for the algebra $\BB$, so the surjection $\U \to \S$ factors through
the algebra $\BB$.  In particular, this gives a surjection $\BB \to
\S$.  Thus, as in the classical case, to complete the proof of Theorem
\ref{qthm:present:S} it is enough to produce a spanning set in $\BB$
of cardinality equal to the dimension of $\S$.  We know from
Proposition \ref{prop:Omega:qspans} that $\qbasis_+$ and $\qbasis_-$
span the algebra $\BB$, and it is clear that these sets have
cardinality equal to the dimension of $\S$, since they are in
bijective correspondence with the set $\basis_+$.  This proves that
$\BB \simeq \S$ and Theorem \ref{qthm:present:S} follows.  It also
follows that the sets $\qbasis_+$ and $\qbasis_-$ are $\Q(v)$-bases
for $\S$.  Hence these sets are linearly independent over $\Q(v)$, and
thus also over $\A$.  Thus they are $\A$-bases for $\BB_\A$.  By
\cite{Du}, the restriction to $\U_\A$ of the map $\U \to \S$ gives a
surjection $\U_\A \to \S_\A$.  It follows that the restriction map
$\BB_\A \to \S_\A$ is an isomorphism (of $\A$-algebras), and that the
sets $\qbasis_+$, $\qbasis_-$ are $\A$-bases for the $q$-Schur algebra
$\S_\A$. Moreover, by Proposition \ref{prop:qidempotent}, $K_i$ and
$K_i^{-1}$ lie in the subalgebra of $\BB_\A$ generated by the
$\sqbinom{K_i}{b}$.  This proves Theorem \ref{qthm:idemp:basis}.

\begin{remark}
We conjecture that for any fixed $i_0$ ($1\le i_0 \le n$) 
either of the sets
$$
\{ E_A K_B F_C \}, \quad \{ F_A K_B E_C \}
$$
($B\in \N^n_{i_0}, A,C \in \N^{\Phi^+}, |A|+|B|+|C|\le d$) is an
$\A$-basis of $\S_\A(n,d)$.  These are a truncated form of
Lusztig's analogue for $\U_\A$ of Kostant's basis for $U_\Z$.  The
conjecture is true when $n=2$; see \cite{DG:quantum}.
\end{remark}

\section{Proof of Theorem \ref{qthm:idemp:present}}
\label{qsec:truncPBW}

To prove Theorem \ref{qthm:idemp:present}, one sets $K_j =
\sum_\lambda v^{\lambda_j}1_\lambda$, $K_j^{-1} = \sum_\lambda
v^{-\lambda_j}1_\lambda$ and verifies that these elements, along with
the $E_i, F_i$, satisfy the relations \eqref{Q1}--\eqref{Q7}.  The
argument is similar to that given in the proof of Theorem
\ref{thm:idemp:present}.  The details are left to the reader.

\section{Applications} \label{sec:apps}

In this section we apply our main results to study some
subalgebras of $\S$. In what follows, we will focus entirely on
the quantum case as the corresponding results for the classical
case are essentially the same.

Recall the triangular decomposition \eqref{BBA:TD}, which we now write
in the form $\S_{\A}=\S_{\A}^+\S_{\A}^0\S_{\A}^-$ in light of the
identification $\S_\A \cong \BB_\A$.  We consider the plus part
$\S_\A^+$.  We give a new proof for the following result
of R.M.\ Green \cite[Prop.\ 2.3]{RGreen}. 

\begin{prop} \label{RGreen's result}
Let $E_A\in \S_{\A}^+$ and suppose the products are taken in the
box order. Then $E_A=0$ if $|A|>d$. Similarly if $F_C\in
\S_{\A}^-$ and the products are taken in the reverse box order,
then $F_C=0$ if $|C|>d$.
\end{prop}

\begin{proof}
Let $E_A\in \S_{\A}^+$ be given in the box order with $|A|>d$. For
each $j=2,\ldots , n$ set
$$
E_{A_j}=\prod_{i=1}^{j-1}X_{j-i,j}^{(A_{j-i,j})},
$$ 
where we write $X_{j-i,j}$ (resp., $A_{j-i,j}$) short for
$X_{\varepsilon_{j-i} - \varepsilon_j}$ (resp., $A_{\varepsilon_{j-i} - \varepsilon_j}$).  
Then $E_A=E_{A_n}\cdots E_{A_2}$. According to
Proposition \ref{qprop:rootXidemp}, if $E_{A_j}1_\lambda\neq 0$, then
$E_{A_j}1_\lambda= 1_{\mu}E_{A_j}$ where $\mu_j=\lambda_j-\sum
_{r=1}^{j-1}A_{j-r,j}$ and $\mu_s=\lambda_s$ if $s>j$. Therefore for
$E_A 1_{\lambda}\neq 0$ it is necessary that $\lambda_j \geq
\sum_{r=1}^{j-1}A_{j-r,j}$ for all $j=2,\ldots, n$, and thus
$$
\sum_{j=2}^n \lambda_j \geq \sum_{j=2}^n\sum_{r=1}^{j-1} A_{j-r,j}
= |A|.
$$ 
This inequality cannot be satisfied since $\Sigma_{j=2}^n\lambda_j
\leq d$ and $|A|>d$. Therefore $E_A1_{\lambda}=0$ for all $\lambda \in
\Lambda(n,d)$.  This forces $E_A=0$, since $1=\sum 1_{\lambda}$. This
completes the proof for the claim about the plus part. The proof for
the other claim is similar.
\end{proof}

\begin{prop}\label{S+basis}
Fix an order on $\Phi^+$. The set of $E_A$ (resp., $F_A$) such that
$|A|\leq d$, with products of factors taken in the designated order,
is an $\A$-basis of $\S_{\A}^+$ (resp., $\S_{\A}^-$).
\end{prop}

\begin{proof}
First consider the box order on $\Phi^+$.  From Lusztig
\cite[\sect4.7]{Lusztig} we know that the set of all $E_A$ (products
taken in the box order) is an $\A$-basis of $\U^+$.  Thus, by
Proposition \ref{RGreen's result}, the set $\Gamma$ consisting of
those $E_A$ satisfying the condition $|A|\leq d$ is an $\A$-spanning
set of $\S_{\A}^+$. The cardinality of this spanning set is equal to
$\dim_{\Q(v)} \S^+$, so it follows that it forms a basis (over
$\Q(v)$) of $\S^+$. Hence the elements in the set are linearly
independent over $\Q(v)$, and thus also over $\A$. Hence they form an
$\A$-basis of $\S_\A^+$.

Now we form products $E_A$ with respect to an {\em arbitrary} (but
fixed) order on $\Phi^+$. By \cite[Theorem 2.4]{Xi2} the set of such
products spans $\S_\A^+$. But the commutation formulas
\eqref{Uq:rel++} do not increase degree; hence when we express an
element of $\Gamma$ as an $\A$-linear combination of $E_A$'s (in the
given fixed order on factors) the degree cannot increase.  It follows
that the set of $E_A$ ($|A|\le d$) must span $\S_\A^+$ (over $\A$).
It follows that this set is also an $\A$-basis of $\S_\A^+$.

This proves the statement for the plus part.  The proof for the minus
part is similar.
\end{proof}

Next we consider the Borel Schur algebras $\S_{\A}^{\geqslant 0} =
\S_\A^0 \S_A^+$ and $\S_{\A}^{\leqslant 0}= \S_\A^0 \S_\A^-$.

\begin{prop}\label{Borelbasis}
Fix an order on $\Phi^+$. With products taken in the specified order,
the set $\{ E_A 1_{\lambda} \mid \ct(E_A) \preceq \lambda \}$ (resp.,
$\{ 1_{\lambda} F_C \mid \ct(F_C) \preceq \lambda \}$) is an
$\A$-basis of $\S_{\A}^{\geqslant 0}$ (resp., $\S_{\A}^{\leqslant 0}$).
\end{prop}

\begin{proof}
From the preceding result and the decomposition $\S_\A^{\geqslant 0} =
\S_\A^+ \S_\A^0$ we see that the set of all $E_A 1_\lambda$ spans
$\S_\A^{\geqslant 0}$.  We can argue as in the proof of Proposition
\ref{prop:Omega:qspans} that with the restriction $\ct(E_A) \preceq
\lambda$ the set still spans.  Clearly this set is linearly
independent over $\A$ since it is a subset of an $\A$-basis of
$\S_\A$.  This proves the first statement.  The proof of the other case
is similar.
\end{proof}

\begin{remark}
We conjecture that the set of all elements of the form $E_A H_B$
($A\in \N^{\Phi^+}$, $B\in \N_{i_0}^n$) satisfying $|A|+|B| \le d$ is
an $\A$-basis of $\S_\A^{\geqslant 0}$, with a similar statement
applying to $\S_\A^{\leqslant 0}$.  
\end{remark}

\bigskip
Now we consider an application to Hecke algebras. Suppose that
$n\ge d$.  Let $\omega = (1^d)$.  Then the subalgebra $1_\omega
\S(n,d) 1_\omega$ is isomorphic with the Hecke algebra
$\mathbf{H}=\mathbf{H}(\Sigma_d)$. If $E_A 1_\lambda F_C$
($\ct(E_AF_C)\preceq \lambda$) is any basis element of $\S$ then by
Propositions \ref{prop:qidempotent} and \ref{qprop:rootXidemp} we see
that $1_\omega E_A 1_\lambda F_C 1_\omega = 0$ unless 
\begin{equation}\label{hecke:condition}
\lambda+
\sum_{\alpha\in \Phi^+} A(\alpha)\alpha = \omega = \lambda+
\sum_{\alpha\in \Phi^+} C(\alpha)\alpha, 
\end{equation}
in which case $1_\omega E_A 1_\lambda F_C 1_\omega = E_A 1_\lambda F_C
= 1_\omega E_AF_C 1_\omega = 1_\omega E_AF_C = E_AF_C 1_\omega$.
Therefore the nonzero elements $E_A 1_\lambda F_C$ of $\qbasis_+$
satisfying condition \eqref{hecke:condition} comprise an $\A$-basis of
$\mathbf{H}$.  There is a similar basis for $\mathbf{H}$ as a subset
of the basis $\qbasis_-$. 

Taking $d=n$, we can see that $\mathbf{H}$ is generated by the
elements $\tgen_i{}= 1_\omega E_i F_i 1_\omega$ ($1 \le i \le n-1$). One
can check directly from relations \eqref{Q1}--\eqref{Q7} and the
propositions in \sect\ref{qsec:B} that these generators satisfy the
following relations:
\begin{align}
&\tgen_i^2 = [2]\, \tgen_i   \tag{H1}\label{H1}\\
&\tgen_i\tgen_j  = \tgen_j\tgen_i  \quad(|i-j|>1)   \tag{H2}\label{H2}\\
&\tgen_i\tgen_{i+1}\tgen_i-\tgen_{i+1}\tgen_{i}\tgen_{i+1} 
= \tgen_i-\tgen_{i+1}. \tag{H3}\label{H3}
\end{align}
Setting $e_i= \tgen_i/[2]$ and putting $q=v^2$ establishes the
equivalence of the above presentation with the presentation in terms
of generators $e_i$ given in Wenzl \cite[\sect2]{Wenzl}.

Note that (with $q=v^2$) the elements $T_i = v^2 - v\tgen_i$ satisfy the
relations
\begin{align}
&T_i^2 = (q-1) T_i + q  \tag{H$1^\prime$}\label{H1'}\\
&T_iT_j  = T_jT_i  \quad(|i-j|>1)  \tag{H$2^\prime$}\label{H2'}\\
&T_i T_{i+1}T_i =
T_{i+1}T_{i}T_{i+1} \tag{H$3^\prime$}\label{H3'}
\end{align}
which is the ``usual'' presentation of $\mathbf{H}$. 

Note also that by \eqref{Q2} and Proposition \ref{qprop:Kmult}(a) we have
$1_\omega E_i F_i 1_\omega = 1_\omega F_i E_i 1_\omega$, so the
alternative ordering of the basis elements of the Schur algebra does
not yield another presentation of $\mathbf{H}$.

\bigskip
\begin{example}
The easiest way to write down the basis elements in the basis
$\qbasis_+$ ($\qbasis_-$ is similar) uses the alternate description in
the $q$-analogue of \eqref{basis:rewrite}, with the idempotent on
the right. Once one has the basis elements, it is a simple matter to
rewrite them with the idempotent anywhere one pleases using
\ref{qprop:rootXidemp} repeatedly.

Given $n,d$ and $\lambda \in \Lambda(n,d)$ set $\qbasis_+(\lambda) = \{
E_A F_C 1_\lambda \mid \ct_R(E_A F_C) \preceq \lambda \}$, so that
$\qbasis_+ = \bigcup_\lambda \qbasis_+(\lambda)$ (disjointly).  The
partition pieces $\qbasis_+(\lambda)$ are obtainable as follows.  One
must choose orders on the factors in $E_AF_C$ (the two parts can be
ordered independently). Once that is done, then for each $j=1\dots,n$
one writes out the set of monomials in variables $X_{ij}$ ($i\ne j$)
of total degree not exceeding $\lambda_j$. (As before, we write
$X_{ij}$ short for $X_\alpha$ with $\alpha = \varepsilon_i -
\varepsilon_j$.) Then one takes the ordered Cartesian product of these
sets over $j$, respecting the given order, with the factor $1_\lambda$
at the right.  (From this it is easy to write a formula for
$|\qbasis_+(\lambda)|$ as a product of binomial coefficients.)

Note that it suffices to describe $\qbasis_+(\lambda)$ just for {\em
dominant} $\lambda \in \Lambda(n,d)$, since the sets indexed by
non-dominant $\lambda$ can be obtained from the dominant one in its
orbit by applying the appropriate permutation to the indices (and then
reordering the product to conform to the specified orders on factors of
$E_A$, $F_C$, if necessary).  We list below the elements $E_AF_C$ such
that $E_A F_C 1_\lambda \in \qbasis_+(\lambda)$, for dominant
$\lambda$.  The elements corresponding with basis elements of the
Hecke algebra are underlined.

For $\S_\A(2,2)$ the elements in question are
\begin{align*}
1_{(2,0)}:\ &\{1, X_{21}, \divided{X_{21}}{2}\}\\
1_{(1,1)}:\ &\{\underline{1}, X_{12}, X_{21}, \underline{X_{12}X_{21}}\}. 
\end{align*}
Thus $\dim \S(2,2) = 2\cdot3 + 4 = 10$. (There are two sets in the
$(2,0)$ orbit, each of cardinality $3$.)

For $\S_\A(3,3)$ we fix the order
$(12)\prec (13)\prec (23) \prec (21) \prec (31) \prec (32)$.  Then the
sets are determined by the elements
\begin{align*}  
1_{(3,0,0)}:\ &\{ 1, X_{21}, X_{31}, \divided{X_{21}}{2}, \divided{X_{31}}{2}, 
 X_{21}X_{31}, \divided{X_{21}}{3}, \divided{X_{31}}{3}, 
\divided{X_{21}}{2}X_{31}, X_{21}\divided{X_{31}}{2} \}\\[3pt]
1_{(2,1,0)}:\ &\{ 1, X_{12}, X_{21}, X_{31}, X_{32}, X_{12}X_{21}, 
X_{12}X_{31}, \divided{X_{21}}{2}, \divided{X_{31}}{2}, X_{21}X_{31},\\
   & X_{21}X_{32}, X_{31}X_{32},  
X_{12}\divided{X_{21}}{2}, X_{12}\divided{X_{31}}{2}, X_{12}X_{21}X_{31}, 
\divided{X_{21}}{2}X_{32},\\ \displaybreak[0]
   & \divided{X_{31}}{2}X_{32}, X_{21}X_{31}X_{32} \}\\[3pt]\displaybreak[0]
1_{(1,1,1)}:\ &\{ \underline{1}, X_{12}, X_{13}, X_{23}, X_{21}, X_{31}, X_{32},
X_{12}X_{13}, X_{12}X_{23}, \underline{X_{12}X_{21}}, \\ \displaybreak[0]
   & X_{12}X_{31}, X_{13}X_{21}, \underline{X_{13}X_{31}}, 
X_{13}X_{32}, X_{23}X_{21}, 
X_{23}X_{31}, \underline{X_{23}X_{32}}, \\ \displaybreak[0]
   & X_{21}X_{32}, X_{31}X_{32}, X_{12}X_{13}X_{21}, X_{12}X_{13}X_{31}, 
X_{12}X_{23}X_{21}, \\
   &\underline{X_{12}X_{23}X_{31}}, \underline{X_{13}X_{21}X_{32}},
X_{13}X_{31}X_{32}, X_{23}X_{21}X_{32}, X_{23}X_{31}X_{32} \}.
\end{align*}
Note that for the $\lambda=(2,1,0)$ case we took the ordered Cartesian
product of $\{1, X_{12}, X_{32} \}$ with $\{1, X_{21}, X_{31},
\divided{X_{21}}{2}, \divided{X_{31}}{2}, X_{21}X_{31} \}$ and for the
$\lambda=(1,1,1)$ case we computed the ordered Cartesian product of
the sets $\{1,X_{13},X_{23}\}$, $\{1, X_{12}, X_{32} \}$, and $\{1,
X_{21}, X_{31} \}$.  There are $10$ elements in the $(3,0,0)$ piece,
$18$ in the $(2,1,0)$ piece, and $27$ in the $(1,1,1)$ piece.  There
are $3$ pieces in the $(3,0,0)$-orbit, $6$ in the $(2,1,0)$-orbit, and
$1$ in the $(1,1,1)$-orbit.  Thus $\dim \S(3,3) =
3\cdot 10 + 6\cdot 18 + 27 = 165$.

\end{example}

\bibliographystyle{amsalpha}

\begin{thebibliography}{100}\frenchspacing\raggedright
\small 

\bibitem[\sf BLM]{BLM} A.A. Beilinson, G. Lusztig, and R. MacPherson,
A geometric setting for the quantum deformation of $\GL_n$, {\em Duke
Math. J.} {\bf61} (1990), 655--677.

\bibitem[\sf CL]{CL} R.W. Carter and G. Lusztig, On the modular
representations of the general linear and symmetric groups, {\em
Math. Z.} {\bf136} (1974), 193--242.

\bibitem[\sf DJ1]{DJ1} R. Dipper and G.D. James, The $q$-Schur algebra,
{\em Proc. London Math. Soc.} {\bf 59} (1989), 23--50.

\bibitem[\sf DJ2]{DJ2} R. Dipper and G.D. James, $q$-tensor space and
$q$-Weyl modules, {\em Trans. Amer. Math. Soc.} {\bf 327} (1991),
251--282.

\bibitem[\sf Do]{Donkin:SA3} S. Donkin, On Schur algebras and related
algebras III: integral representations, {\em Math. Proc. Cambridge
Philos. Soc.} {\bf116} (1994), 37--55.

\bibitem[\sf DG1]{DG} S. Doty and A. Giaquinto, Presenting Schur
algebras as quotients of the universal enveloping algebra of $\gl_2$,
{\em Algebras and Representation Theory}, to appear. 

\bibitem[\sf DG2]{DG:quantum} S. Doty and A. Giaquinto, Presenting
quantum Schur algebras as quotients of the quantized enveloping
algebra of $\gl_2$, preprint, Loyola University Chicago, December
2000.

\bibitem[\sf DG3]{DG:announce} S. Doty and A. Giaquinto, Generators
and relations for Schur algebras, {\em Electronic Research
Announc. Amer. Math. Soc.} {\bf7} (2001), 54--62.


\bibitem[\sf Du]{Du} Jie Du, A note on quantized Weyl reciprocity at
roots of unity, {\em Algebra Colloq.} {\bf2} (1995), 363--372.

\bibitem[\sf Gr]{Green} J.~A.~Green, {\em Polynomial Representations
of $\GL_n$}, (Lecture Notes in Math.~{\bf830}), Springer-Verlag, New
York 1980.

\bibitem[\sf RG1]{RG:thesis} R.M. Green, Ph.D. thesis, University of
Warwick, 1995.

\bibitem[\sf RG2]{RGreen} R.M. Green, $q$-Schur algebras as quotients
of quantized enveloping algebras, {\em J. Algebra} {\bf185} (1996),
660--687.

\bibitem[\sf Ji]{Jimbo} M. Jimbo, A $q$-analogue of $U(\gl(N+1))$,
Hecke algebra, and the Yang-Baxter equation, {\em Letters
Math. Physics} {\bf11} (1986), 247--252.

\bibitem[\sf Ko]{Kostant} B. Kostant, Groups over $\Z$, {\em
Proc. Symposia Pure Math.} {\bf9} (1966), 90--98.

\bibitem[\sf Lu1]{Lusztig:88} G. Lusztig, Quantum deformations of
certain simple modules over enveloping algebras, {\em Advances in
Math.} {\bf70} (1988), 237--249.

\bibitem[\sf Lu2]{Lusztig} G. Lusztig, Finite dimensional Hopf
algebras arising from quantized universal enveloping algebras,
{\em J. Amer. Math. Soc.} {\bf 3} (1990), no. 1, 257--296.

\bibitem[\sf Lu3]{Lusztig:book} G. Lusztig, {\em Introduction to
Quantum Groups}, Birkh\"{a}user Boston 1993.

\bibitem[\sf Ro]{Rosso} M. Rosso,
An analogue of P.B.W. theorem and the universal $R$-matrix for
$U\sb h{\rm sl}(N+1)$. {\em Comm. Math. Phys.} {\bf 124} (1989),
no. 2, 307--318.

\bibitem[\sf Sc]{Schur} I. Schur, \"Uber die rationalen Darstellungen
der allgemeinen linearen Gruppe, 1927; reprinted in: I. Schur, {\em
Gesammelte Abhandlungen}, Vol. III, pp. 68--85, Springer-Verlag,
Berlin, 1973.

\bibitem[\sf We]{Wenzl} H. Wenzl, Hecke algebras of type $A_n$ and
subfactors, {\em Invent. Math.} {\bf92} (1988), 349--383.

\bibitem[\sf Xi1]{Xi} Nanhua Xi, Root vectors in quantum groups,
{\em Comment. Math. Helv.} {\bf69} (1994), 612--639.

\bibitem[\sf Xi2]{Xi2} Nanhua Xi,  A commutation formula for root
vectors in quantized enveloping algebras, {\em  Pacific J. Math.}
{\bf 189} (1999), no. 1, 179--199.

\end{thebibliography}

\end{document}